\documentclass[english]{elsarticle}
\usepackage[T1]{fontenc}
\usepackage[latin9]{inputenc}
\usepackage{float}
\usepackage{enumitem}
\usepackage{amsmath}
\usepackage{amsthm}
\usepackage{amssymb}
\usepackage{graphicx}

\makeatletter

\providecommand{\tabularnewline}{\\}

\newlength{\lyxlabelwidth}      
  \theoremstyle{definition}
  \newtheorem{defn}{\protect\definitionname}
\newenvironment{elabeling}[2][]%
{\settowidth{\lyxlabelwidth}{#2}
\begin{description}[font=\normalfont,style=sameline,
leftmargin=\lyxlabelwidth,#1]}
{\end{description}}

\@ifundefined{showcaptionsetup}{}{%
 \PassOptionsToPackage{caption=false}{subfig}}
\usepackage{subfig}
\makeatother

\usepackage{babel}
  \providecommand{\definitionname}{Definition}

\begin{document}

\begin{frontmatter}{}

\title{Exact Nonlinear Model Reduction for a von Kármán beam: \\Slow-Fast
Decomposition and Spectral Submanifolds}

\author{Shobhit Jain, Paolo Tiso and George Haller\footnote{Corresponding author. Email: georgehaller@ethz.ch}}

\address{Institute for Mechanical Systems, ETH Zürich,}

\address{Leonhardstrasse 21, 8092 Zürich, Switzerland}
\begin{abstract}
We apply two recently formulated mathematical techniques, Slow-Fast
Decomposition (SFD) and Spectral Submanifold (SSM) reduction, to a
von Kármán beam with geometric nonlinearities and viscoelastic damping.
SFD identifies a global slow manifold in the full system which attracts
solutions at rates faster than typical rates within the manifold.
An SSM, the smoothest nonlinear continuation of a linear modal subspace,
is then used to further reduce the beam equations within the slow
manifold. This two-stage, mathematically exact procedure results in
a drastic reduction of the finite-element beam model to a one-degree-of
freedom nonlinear oscillator. We also introduce the technique of spectral
quotient analysis, which gives the number of modes relevant for reduction
as output rather than input to the reduction process. 
\end{abstract}

\end{frontmatter}{}

\section{Introduction}

Computer simulations are routinely performed in today's technological
world for modeling and response prediction of almost any physical
phenomenon. The ever-increasing demand for realistic simulations leads
to a higher level of detail during the modeling phase, which in turn
increases the complexity of the models and results in a bigger problem
size. Typically, such physical processes are mathematically modeled
using partial differential equations (PDEs), which are discretized
(e.g. using Finite Elements, Finite differences, Finite volumes methods
etc.) to obtain problems with a finite (but usually large) number
of unknowns. Despite the tremendous increase in computational power
over the past decades, however, the time required to solve high-dimensional
discretized models remains a bottleneck towards efficient and optimal
design of structures. Model order reduction (MOR) aims to reduce the
computational efforts in solving such large problems.

The classical approach to model reduction involves a linear projection
of the full system onto a set of basis vectors. This linear projection
is characterized by a matrix whose columns span a suitable low-dimensional
subspace. Various techniques have been applied to high-dimenisional
systems to obain such a reduction basis, including the Proper Orthogonal
Decomposition (POD) \cite{POD1,POD2,POD3} (also knowns as Singular
Value Decomposition (SVD), Karhunen-Loeve Decomposition), Linear Normal
Modes (LNM) and Krylov subspace projection \cite{krylov}. Once a
suitable basis is chosen, the reduced-order model (ROM) is then obtained
using Galerkin projection. Similar linear projection techniques have
been devised for component-mode synthesis (CMS), such as the Craig-Bampton
method \cite{Craig Bampton}. An implicit assumption to all linear
projection techniques is that the full system dynamics evolves in
a lower-dimensional linear invariant \textit{subspace} of the phase
space of the system. While such linear subspaces do exist for linear
systems (linear modal subspaces), they are generally non-existent
in nonlinear systems. This results in a priori unknown and potentially
large errors for linear projection-based reduction methods, necessitating
the verification of the accuract of the reduction procedure on a case-by-case
basis. 

More recent trends in model reduction account for this issue by constructing
the reduced solution over nonlinear manifolds \cite{qm,qm2,NLDR}.
The seminal idea of Shaw and Pierre \cite{shawpierre93} is to construct
assumed nonlinear invariant surfaces (nonlinear normal modes) that
act as continuations of linear normal mode families under the addition
of nonlinear terms near an equilibrium (see \cite{kerschen2014} for
a recent review). More heuristic reduction procedures in structural
dynamics include the static condensation approaches, where the fast
or stiff variables in the system are intuitively identified, and statically
enslaved to the slow or flexible ones (cf. \cite{mignolet} for a
review). Guyan reduction \cite{Guyan} is a classic example of this
reduction philosophy at the linear level. 

Most classic reduction techniques, coming from a intuive and heuristic
standpoint, do not provide an a priori estimate of their accuracy
or even validity for a given system. To this effect, Haller and Ponsioen
\cite{haller16sf} proposed requirements for mathematically justifiable
and robust model reduction in a non-linear mechanical system. These
requirements ensure not only that the lower dimensional attracting
structure (manifold) in the phase space is robust under perturbation,
but also that the full system trajectories are attracted to it and
synchronize with the reduced model trajectories at rates that are
faster than typical rates within the manifold. 

Recent advances in nonlinear dynamics enable such an \textit{exact}
model-reduction using the slow-fast decomposition (SFD) \cite{haller16sf}
and spectral submanifold (SSM) based reduction \cite{haller16}. SFD
is a general procedure to identify if a mechanical system exhibits
a global partitioning of degrees of freedom into slow (flexible) and
fast (stiff) components such that the fast variables can be enslaved
to the slow ones. This results in a global ROM containing only the
slow degrees of freedom of the full system. SSMs, on the other hand,
are the smoothest nonlinear extensions of the linear modal subspaces
near an asymptotically stable equilibrium. Neither SFD nor SSM-based
reduction has been applied to problems with high numbers of degrees
of freedom. A beam model is often the first step in showing the potential
of a new technique for reduction of high dimensional systems (cf.
\cite{ilbeigi16,craig89,nnmbeam}). To this effect, we combine here,
for the first time, the application of these techniques on a finite-element
discretized nonlinear von Kármán beam model.

We first show that the beam model satisfies the requirements of SFD.
The corresponding ROM is subsequently obtained on a slow manifold
defined over the transverse degrees of freedom of the beam. This SFD-reduced
ROM posseses no clear spectral gaps but a spectral quotient analysis
nonetheless enables a further reduction to an SSM using the formulas
given by Szalai et. al \cite{Szalai16}. Importantly, our spectral
quotient analysis returns the number of modes relevant for reduction,
instead of postulating or deducing this number from numerical experimentation.
In the end, our two-step reduction results in an exact ROM with a
drastic reduction in the number of degrees of freedom of the system. 

In the next section, we start by reviewing the main steps involved
in the derivation of the governing PDEs for the von-Karman beam and
the non-dimensionalization we performed upon them. The finite-element
disretized equations obtained from the resulting nondimensionalized
PDEs are then presented. This system of equations is then first reduced
using the SFD in Section 3. The SSM-based reduction applied to the
SFD-reduced system in Section 4. The conclusions, along with scope
for further work, are presented in Section 5. 
\begin{center}
\begin{table}[h]
\begin{centering}
\begin{tabular}{|c|c|}
\hline 
Symbol & Meaning (unit)\tabularnewline
\hline 
\hline 
$L$ & Length of beam (m)\tabularnewline
\hline 
$h$ & Height of beam (m)\tabularnewline
\hline 
$b$ & Width of beam (m)\tabularnewline
\hline 
$A$ & Area of cross section (m$^{2}$) = $bh$\tabularnewline
\hline 
$E$ & Young's Modulus (Pa)\tabularnewline
\hline 
$\kappa$ & Viscous damping rate of material (Pa s)\tabularnewline
\hline 
$\rho$ & Density (kg/m$^{3}$)\tabularnewline
\hline 
$\tau$ & Non-dimentionalized time\tabularnewline
\hline 
\end{tabular}
\par\end{centering}
\caption{Notation}
\end{table}
\par\end{center}

\section{Setup}

\begin{figure}[h]
\begin{centering}
\includegraphics[width=0.8\textwidth]{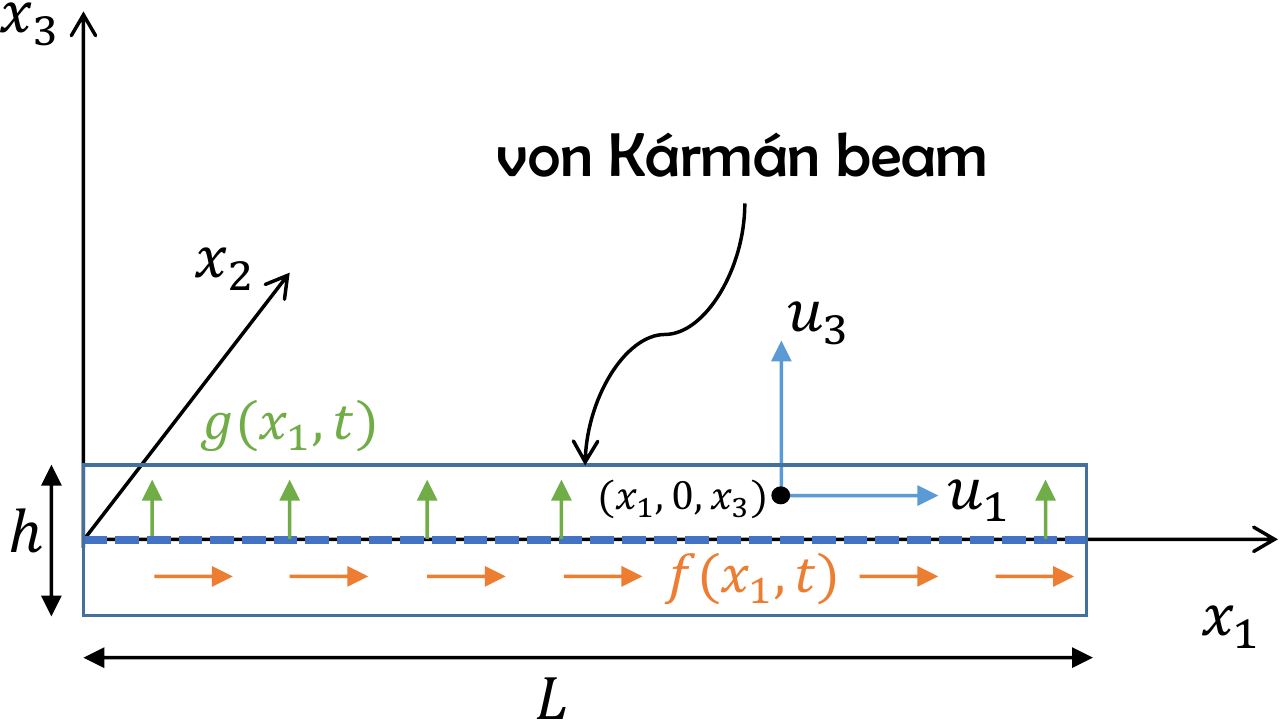}
\par\end{centering}
\caption{\label{fig:beam}The schematic of a two-dimensional von Kármán beam
with height $h$ and length $L$, initially aligned with the $x_{1}$.
The $x_{1}$ and $x_{3}$ displacements of a material point with coordinates
$(x_{1},0,x_{3})$ are given by $u_{1}$ and $u_{3}$, respectively.
The transverse and axial load per unit length applied to the beam
are given by $g(x_{1},t)$ and $f(x_{1},t)$, respectively.}
\end{figure}

We briefly summarize the main steps leading to the derivation of the
governing partial differential equations (PDEs) for the von Kármán
beam (see, e.g., Reddy \cite{reddy10} for a detailed derivation).
We consider a straight 2D beam aligned initially with the $x_{1}$
axis, as shown in Figure \ref{fig:beam}. The motion of the beam takes
place in the $x_{1}-x_{3}$ plane. 
Assuming the Euler-Bernoulli hypothesis for the kinematics of bending
(i.e., that the lines initially straight and perpendicular to the
beam axis remain so after deformation), we obtain the displacement
field: 

\begin{align}
u_{1}(x_{1},x_{3}) & =u_{0}(x_{1})-x_{3}\partial_{x_{1}}w_{0}(x_{1}),\nonumber \\
u_{2}(x_{1},x_{3}) & =0,\nonumber \\
u_{3}(x_{1},x_{3}) & =w_{0}\left(x_{1}\right),\label{eq:Kinematics}
\end{align}
where $\left(u_{1},u_{2},u_{3}\right)$ denote the $\left(x_{1},\,x_{2},\,x_{3}\right)$
displacements of a material point with coordinates $(x_{1},0,x_{3})$,
and $\left(u_{0}(x_{1}),\,w_{0}(x_{1})\right)$ are the $\left(x_{1},\,x_{3}\right)$
displacements of any material point lying on the reference line given
by the $x_{1}$ axis. The von Kármán strain approximation of the Green-Lagrange
strain for moderate rotations is given by
\begin{align}
\varepsilon_{11} & =\partial_{x_{1}}u_{0}+\frac{1}{2}(\partial_{x_{1}}w_{0})^{2}-z\partial_{x_{1}}^{2}w_{0}.\label{eq:VKstrain}
\end{align}
Using the virtual work principle, we can formulate the equations of
motion in terms of the primary unknowns $u_{0,\,}w_{0}$ as
\begin{align*}
\rho A\mathrm{\partial}_{t}^{2}u_{0}= & \,\partial_{x_{1}}\underbrace{\left[\int_{A}\sigma~\mathrm{d}A\right]}_{N}+f(x_{1},t)\,,\\
\rho A\mathrm{\partial}_{t}^{2}w_{0}= & \,\partial_{x_{1}}^{2}\left[\int_{A}x_{3}\sigma~\mathrm{d}A\right]+\partial_{x_{1}}\left(N\partial_{x_{1}}w_{0}\right)+g(x_{1},t)\,,
\end{align*}
where $f$ and $g$ are external body forces per unit length of the
beam in the $x_{1}$ and $x_{3}$ directions, respectively; $A$ is
the area of the cross-section; and $\sigma$ is the $\sigma_{11}$
component of the Cauchy stress. We choose the Kelvin-Voigt model for
viscoelasticity as constitutive law to relate the stress $\sigma$
to the von Kármán strain $\varepsilon_{11}$(\ref{eq:VKstrain}) and
to the corresponding strain rate as
\[
\sigma=E\varepsilon+\kappa\dot{\varepsilon}\,.
\]
Here $E$ denotes Young's modulus and $\kappa$ is the rate of viscous
damping for the material. Further assumptions include a uniform rectangular
cross section with the reference $x_{1}$ axis passing through the
centroid of the cross-section, and a non-dimensionalization of variables
as $x=\frac{x_{1}}{L},\,w=\frac{w_{0}}{h},\,u=\frac{u_{0}}{h},\:\tau=\frac{h}{L^{2}}\sqrt{\frac{E}{\rho}}t,\,p(x,\tau)=\frac{\beta L^{4}}{bEh^{4}}f\left(x,\frac{L^{2}}{h}\sqrt{\frac{\rho}{E}}\tau\right),\,q(x,\tau)=\frac{\alpha L^{4}}{bEh^{4}}g\left(x,\frac{L^{2}}{h}\sqrt{\frac{\rho}{E}}\tau\right)$
($\alpha,\beta$ have been introduced as scaling factors for transverse
and axial loading respectively). These lead to the following dimensionless
PDEs governing the beam behavior:
\begin{gather}
\ddot{w}+\frac{1}{12}\partial_{x}^{4}w+\frac{\zeta\epsilon}{12}\partial_{x}^{4}\dot{w}-\frac{1}{\epsilon}\partial_{x}\left(\partial_{x}u\,\partial_{x}w\right)-\zeta\partial_{x}(\partial_{x}\dot{u\,}\partial_{x}w)\nonumber \\
-\frac{1}{2}\partial_{x}\left(\partial_{x}w\right)^{3}-\zeta\epsilon\partial_{x}\left(\left(\partial_{x}w\right)^{2}\partial_{x}\dot{w}\right)=\alpha q(x,\tau)\,,\label{eq:govPDE}\\
\ddot{u}-\frac{1}{\epsilon}\partial_{x}^{2}u-\frac{\zeta}{\epsilon}\partial_{x}^{2}\dot{u}-\frac{1}{2\epsilon}\partial_{x}\left(\partial_{x}w\right)^{2}-\zeta\partial_{x}\left(\partial_{x}w\partial_{x}\dot{w}\right)=\beta p(x,\tau)\,.\nonumber 
\end{gather}
Here $\dot{(\bullet)}$ denotes $\partial_{\tau}(\bullet)$, $\epsilon=\frac{h}{L}$
is the length to thickness ratio and $\zeta=\frac{\kappa\rho^{1/2}}{E^{3/2}L}$
is a dimensionless constant resulting from non-dimensionalization.
We have chosen the displacements to be scaled with respect to the
thickness $h$ of the beam, and not the length. Moreover, the applied
loading has been nondimensionalized and scaled with respect to a load
which leads to transverse displacements in the order of the thickness
$h$ of the beam. This is because the von Kármán kinematic approximations
are justified for displacements and forces in this range. Furthermore,
the time $t$ has been nondimensionalized with respect to a time period
that is representatitve of natural frequency of the oscillation of
the beam, since the structural response is generally studied at such
time scales. 

Upon finite-element discretization of the non-dimensional system (\ref{eq:govPDE})
with cubic shape functions for $w$ and linear shape functions for
$u$ (see, e.g., Crisfield \cite{crisfield96}), we obtain the finite-dimensional
discretized version of (\ref{eq:govPDE}) as
\begin{gather}
\mathbf{\mathbf{M}}_{1}\ddot{\mathbf{x}}+\zeta\epsilon\left(\mathbf{K}_{1}+\boldsymbol{\mathcal{C}}(\mathbf{x})\right)\dot{\mathbf{x}}+\zeta\boldsymbol{\mathcal{D}}(\mathbf{x})\dot{\mathbf{y}}+\mathbf{K}_{1}\mathbf{x}+\dfrac{1}{\epsilon}\boldsymbol{\mathcal{F}}(\mathbf{x},\mathbf{y})+\boldsymbol{\mathcal{G}}(\mathbf{x})=\mathbf{\alpha q}(\tau)\,,\nonumber \\
\mathbf{\mathbf{M}}_{2}\ddot{\mathbf{y}}+\frac{\zeta}{\epsilon}\mathbf{K}_{2}\dot{\mathbf{y}}+\zeta\boldsymbol{\mathcal{E}}(\mathbf{x})\dot{\mathbf{x}}+\frac{1}{\epsilon^{2}}\mathbf{K}_{2}\mathbf{y}+\dfrac{1}{\epsilon}\boldsymbol{\mathcal{H}}(\mathbf{x})=\beta\mathbf{p}(\tau),\label{eq:FEdisc}
\end{gather}
where $\mathbf{x}\in\mathbb{R}^{n_{s}},\mathbf{y}\in\mathbb{R}^{n_{f}}$
are the finite dimensional (discretized) counterparts of the unknowns
$w,u$ respectively ($n_{s},\,n_{f}$ being the number of unknowns
dependent on the finite-element discretization), and $\mathbf{\mathbf{M}}_{1},\mathbf{\mathbf{K}}_{1}\in\mathbb{R}^{n_{s}\times n_{s}}$
and $\mathbf{\mathbf{M}}_{2},\mathbf{\mathbf{K}}_{2}\in\mathbb{R}^{n_{f}\times n_{f}}$
are the corresponding mass and stiffness matrices. Here, $\boldsymbol{\mathcal{F}}$
(a bilinear function of the form $\left(\boldsymbol{\mathcal{F}}(\mathbf{x},\mathbf{y})\right)_{i}=F_{ijk}x_{j}y_{k},\,i,j\in\{1,\dots,n_{s}\},\,k\in\{1,\dots,n_{f}\}$),
$\boldsymbol{\mathcal{G}}$ (a cubic function of the form $\left(\boldsymbol{\mathcal{G}}(\mathbf{x})\right)_{i}=G_{ijkl}x_{j}x_{k}x_{l},\,i,j,k,l\in\{1,\dots,n_{s}\}$)
and $\boldsymbol{\mathcal{H}}$ (a quadratic function of the form
$\left(\boldsymbol{\mathcal{H}}(\mathbf{x})\right)_{i}=H_{ijk}x_{j}x_{k},\,i\in\{1,\dots,n_{f}\},\,j,k\in\{1,\dots,n_{s}\}$)
correspond to the nonlinear elastic forces in the beam; $\boldsymbol{\mathcal{C}}$
(a quadratic function of the form $\left(\boldsymbol{\mathcal{C}}(\mathbf{x})\right)_{ij}=C_{ijkl}x_{k}x_{l},\,i,j,k,l\in\{1,\dots,n_{s}\}$),
$\boldsymbol{\mathcal{D}}$ (a linear function of the form $\left(\boldsymbol{\mathcal{D}}(\mathbf{x})\right)_{ij}=D_{ijk}x_{k},\,i,k\in\{1,\dots,n_{s}\},\,j\in\{1,\dots,n_{f}\}$)
and $\boldsymbol{\mathcal{E}}$ (a linear function of the form $\left(\boldsymbol{\mathcal{E}}(\mathbf{x})\right)_{ij}=E_{ijk}x_{k},\,j,k\in\{1,\dots,n_{s}\},\,i\in\{1,\dots,n_{f}\}$)
correspond to the nonlinear contributions resulting from the viscoelastic
material damping.

\section{Slow-Fast Decomposition}

\subsection{Verification of assumptions for the application of SFD}

In the finite-element discretized system (\ref{eq:FEdisc}), the $\mathbf{y}$
variables (representing the axial displacement components) are stiffer
and hence faster than the $\mathbf{x}$ variables (representing the
transverse displacement components). Such a global difference of speeds
indicates the possible existence of a lower-dimensional slow manifold,
as described in mathematical terms by the geometric singular perturbation
theory of Fenichel \cite{fenichel79}. If such a slow manifold is
robust and attracts nearby solutions, then the dynamics on this manifold
provides and exact reduced-order model with which all nearby solutions
synchronize exponentially fast.

For general finite-dimensional mechanical systems characterized by
such a dicotomy of time scales, Haller and Ponsioen \cite{haller16sf}
deduced conditions under which positions and velocities in the fast
degrees of freedom ($\mathbf{y}$,$\dot{\mathbf{y}}$) can be expressed
as a graph over their slow counterparts ($\mathbf{x}$,$\dot{\mathbf{x}}$),
resulting in a globally exact model reduction. If these conditions
for a Slow-Fast Decomposition (SFD) are satisfied, then all trajectories
of the full system (close enough to the slow manifold in the phase
space) synchronize with the reduced model trajectories at rates faster
than those within the slow manifold. To check these conditions, we
take $0<\epsilon\ll1$, the thickness to length ratio of the beam,
as the required non-dimensional small parameter. The system (\ref{eq:FEdisc})
can then be reformulated in terms of $\frac{\mathbf{y}}{\epsilon}$
as

\begin{gather}
\mathbf{\mathbf{M}}_{1}\ddot{\mathbf{x}}+\zeta\epsilon\left(\mathbf{K}_{1}+\boldsymbol{\mathcal{C}}(\mathbf{x})\right)\dot{\mathbf{x}}+\zeta\boldsymbol{\mathcal{D}}(\mathbf{x})\dot{\mathbf{y}}+\mathbf{K}_{1}\mathbf{x}+\boldsymbol{\mathcal{F}}\left(\mathbf{x},\frac{\mathbf{y}}{\epsilon}\right)+\boldsymbol{\mathcal{G}}(\mathbf{x})=\mathbf{\alpha q}(\tau)\,,\nonumber \\
\mathbf{\mathbf{M}}_{2}\ddot{\mathbf{y}}+\frac{\zeta}{\epsilon}\mathbf{K}_{2}\dot{\mathbf{y}}+\zeta\boldsymbol{\mathcal{E}}(\mathbf{x})\dot{\mathbf{x}}+\frac{1}{\epsilon}\mathbf{K}_{2}\frac{\mathbf{y}}{\epsilon}+\dfrac{1}{\epsilon}\boldsymbol{\mathcal{H}}(\mathbf{x})=\mathbf{\beta p}(\tau).\label{eq:SFD-1}
\end{gather}
The mass-normalized forcing terms are defined in terms of the new
variable $\boldsymbol{\eta}=\frac{\mathbf{y}}{\epsilon}$ as 
\begin{align}
\mathbf{P}_{1}(\mathbf{x},\dot{\mathbf{x}},\boldsymbol{\eta},\dot{\mathbf{y}},\tau;\epsilon) & =-\mathbf{\mathbf{M}}_{1}^{-1}\left(\zeta\epsilon\left(\mathbf{K}_{1}+\boldsymbol{\mathcal{C}}(\mathbf{x})\right)\dot{\mathbf{x}}+\zeta\boldsymbol{\mathcal{D}}(\mathbf{x})\dot{\mathbf{y}}+\mathbf{K}_{1}\mathbf{x}+\mathcal{\boldsymbol{F}}\left(\mathbf{x},\boldsymbol{\eta}\right)+\boldsymbol{\mathcal{G}}(\mathbf{x})-\mathbf{\alpha q}(\tau)\right)\,,\nonumber \\
\mathbf{P}_{2}(\mathbf{x},\dot{\mathbf{x}},\boldsymbol{\eta},\dot{\mathbf{y}},\tau;\epsilon) & =-\epsilon\mathbf{\mathbf{M}}_{2}^{-1}\left(\frac{\zeta}{\epsilon}\mathbf{K}_{2}\dot{\mathbf{y}}+\zeta\boldsymbol{\mathcal{E}}(\mathbf{x})\dot{\mathbf{x}}+\frac{1}{\epsilon}\mathbf{K}_{2}\boldsymbol{\eta}+\dfrac{1}{\epsilon}\boldsymbol{\mathcal{H}}(\mathbf{x})-\mathbf{\beta p}(\tau)\right)\label{eq:SFD-2}\\
 & =-\mathbf{M}_{2}^{-1}\left(\zeta\mathbf{K}_{2}\dot{\mathbf{y}}+\epsilon\zeta\boldsymbol{\mathcal{E}}(\mathbf{x})\dot{\mathbf{x}}+\mathbf{K}_{2}\boldsymbol{\eta}+\boldsymbol{\mathcal{H}}(\mathbf{x})-\epsilon\mathbf{\beta p}(\tau)\right).\nonumber 
\end{align}
The main conditions of the SFD procedure, as derived by Haller and
Ponsioen \cite{haller16sf}, are the following.
\begin{description}
\item [{(A1)}] \textbf{Nonsingular extension to $\epsilon=0$:} $\mathbf{P}_{1}$
and $\mathbf{P}_{2}$ should possess smooth (in fact $C^{\infty}$)
extension to their respective $\epsilon=0$ limits, which is the case
here by (\ref{eq:SFD-2}).
\item [{(A2)}] \textbf{Existence of a critical manifold:} The algebraic
equation $\mathbf{P}_{2}(\mathbf{x},\dot{\mathbf{x}},\boldsymbol{\eta},\mathbf{0},\tau;0)\equiv\mathbf{0}$
should be solvable for $\boldsymbol{\eta}$ on an open, bounded domain.
Indeed, for any $\mathcal{D}_{0}\subset\mathbb{R}^{n_{s}}\times\mathbb{R}^{n_{s}}\times\mathcal{T}$
open and bounded, the set $\mathcal{M}_{0}(\tau)$ defined by 
\[
\boldsymbol{\eta}=\mathbf{G}_{0}(\mathbf{x},\dot{\mathbf{x}},\tau):=-\mathbf{K}_{2}^{-1}\boldsymbol{\mathcal{H}}(\mathbf{x})
\]
 satisfies $\mathbf{P}_{2}(\mathbf{x},\dot{\mathbf{x}},\mathbf{G}_{0}(\mathbf{x},\dot{\mathbf{x}},\tau),\mathbf{0},\tau;0)\equiv\mathbf{0}$
for all $(\mathbf{x},\dot{\mathbf{x}},\tau)\in\mathcal{D}_{0}$, forming
a critical manifold $\mathcal{M}_{0}$ in the language of singular
perturbation theory. Clearly, the critical manifold is independent
of $\tau$ in our current setting. If, however, we assumed that$\beta=\mathcal{O}(\frac{1}{\epsilon})$
(i.e., the beam is excited in the axial direction by forces that are
an order of magnitude larger than those in the transverse direction),
then we would obtain a time-dependent critical manifold.
\item [{Notation:}] For clarity, we denote any expression $(\bullet)$
evaluated over the critical manifold $\mathcal{M}_{0}$ by $\ensuremath{\overline{(\bullet)}:=(\bullet)|_{\boldsymbol{\eta}=\mathbf{G}_{0}(\mathbf{x},\dot{\mathbf{x}},\tau),\,\dot{\mathbf{y}}=\mathbf{0},\,\epsilon=0}}$
\item [{(A3)}] \textbf{Asymptotic stability of the critical manifold}:
With the matrices 
\begin{align}
\mathbf{A}(\mathbf{x},\dot{\mathbf{x}},\tau) & =-\partial_{\dot{\mathbf{y}}}\mathbf{P}_{2}(\mathbf{x},\dot{\mathbf{x}},\mathbf{G}_{0}(\mathbf{x},\dot{\mathbf{x}},\tau),\mathbf{0},\tau;0)=-\overline{\partial_{\dot{\mathbf{y}}}\mathbf{P}_{2}}=\zeta\mathbf{M}_{2}^{-1}\mathbf{K}_{2},\nonumber \\
\mathbf{B}(\mathbf{x},\dot{\mathbf{x}},\tau) & =-\partial_{\boldsymbol{\eta}}\mathbf{P}_{2}(\mathbf{x},\dot{\mathbf{x}},\mathbf{G}_{0}(\mathbf{x},\dot{\mathbf{x}},\tau),\mathbf{0},\tau;0)=-\overline{\partial_{\boldsymbol{\eta}}\mathbf{P}_{2}}=\mathbf{M}_{2}^{-1}\mathbf{K}_{2},\label{eq:ABdef}
\end{align}
the equilibrium solution $\boldsymbol{\eta}\equiv\mathbf{0}\in\mathbb{R}^{n_{f}}$
of the unforced, constant-coefficient linear system
\begin{align}
\boldsymbol{\eta}{}^{\prime\prime}+\mathbf{A}(\mathbf{x},\dot{\mathbf{x}},\tau)\boldsymbol{\eta}^{\prime}+\mathbf{B}(\mathbf{x},\dot{\mathbf{x}},\tau)\boldsymbol{\eta}\Leftrightarrow\mathbf{M}_{2}\boldsymbol{\eta}{}^{\prime\prime}+\zeta\mathbf{K}_{2}\boldsymbol{\eta}^{\prime}+\mathbf{K}_{2}\boldsymbol{\eta} & =\mathbf{0}=\mathbf{0}\label{eq:associated linear system}
\end{align}
should be asymptotically stable for all fixed parameter values $(\mathbf{x},\dot{\mathbf{x}},\tau)\in\mathcal{D}_{0}$.
This is again satisfied in our setting since $\mathbf{M}_{2},\mathbf{K}_{2}$
are positive definite matrices and $\zeta>0$.
\end{description}

\subsection{Global reduced-order model from SFD}

As shown by Haller and Ponsioen \cite{haller16sf}, assumptions (A1)-(A3)
guarantee that the critical manifold $\mathcal{M}_{0}(\tau)$ perturbs
into a nearby attracting slow manifold $\mathcal{M}_{\epsilon}(\tau)$
for $\epsilon>0$ small enough. On this slow manifold, the discretized
beam system (\ref{eq:FEdisc}) admits an exact reduced order model
given by
\begin{equation}
\ddot{\mathbf{x}}-\overline{\mathbf{P}_{1}}-\epsilon\left[\overline{\partial_{\boldsymbol{\eta}}\mathbf{P}_{1}}\mathbf{G}_{1}(\mathbf{x},\dot{\mathbf{x}},\tau)+\overline{\partial_{\dot{\mathbf{y}}}\mathbf{P}_{1}}\mathbf{H_{0}}(\mathbf{x},\dot{\mathbf{x}},\tau)+\overline{\partial_{\epsilon}\mathbf{P}_{1}}\right]+\mathcal{O}(\epsilon^{2})=\mathbf{0},\label{eq:SFD-ROM}
\end{equation}
where 
\begin{align*}
\mathbf{H}_{0}(\mathbf{x},\dot{\mathbf{x}},\tau) & =\left[\partial_{\mathbf{x}}\mathbf{G}_{0}(\mathbf{x},\dot{\mathbf{x}},\tau)\right]\dot{\mathbf{x}}+\left[\partial_{\dot{\mathbf{x}}}\mathbf{G}_{0}(\mathbf{x},\dot{\mathbf{x}},\tau)\right]\overline{\mathbf{P}_{1}}+\partial_{\tau}\mathbf{G}_{0}(\mathbf{x},\dot{\mathbf{x}},\tau)\,,\\
\mathbf{G}_{1}(\mathbf{x},\dot{\mathbf{x}},\tau) & =\left[\overline{D_{\boldsymbol{\eta}}\mathbf{P}_{2}}\right]^{-1}\overline{D_{\dot{\mathbf{y}}}\mathbf{P}_{2}}\mathbf{H}_{0}(\mathbf{x},\dot{\mathbf{x}},\tau)
\end{align*}
constitute the higher-order terms in the equations describing the
slow manifold $\mathcal{M}_{\epsilon}(\tau)$ as
\begin{align}
\mathbf{y} & =\epsilon\mathbf{G}_{0}(\mathbf{x},\dot{\mathbf{x}},\tau)+\epsilon^{2}\mathbf{G}_{1}(\mathbf{x},\dot{\mathbf{x}},\tau)+\mathcal{O}\left(\epsilon^{3}\right)\,,\nonumber \\
\dot{\mathbf{y}} & =\epsilon\mathbf{H}_{0}(\mathbf{x},\dot{\mathbf{x}},\tau)+\mathcal{O}\left(\epsilon^{2}\right)\,.\label{eq:SFDenslavement}
\end{align}
 In the context of the beam example considered here, the exact reduced-order
model in (\ref{eq:SFD-ROM}) can be explicitly written out in the
following form:
\begin{gather}
\mathbf{M}_{1}\ddot{\mathbf{x}}+\mathbf{K}_{1}\mathbf{x}+\boldsymbol{\mathcal{F}}\left(\mathbf{x},\mathbf{G}_{0}(\mathbf{x},\dot{\mathbf{x}},\tau)\right)+\boldsymbol{\mathcal{G}}(\mathbf{x})+\nonumber \\
\epsilon\left[\underbrace{\overline{\partial_{\boldsymbol{\eta}}\boldsymbol{\mathcal{F}}\left(\mathbf{x},\boldsymbol{\eta}\right)}\mathbf{G}_{1}(\mathbf{x},\dot{\mathbf{x}},\tau)}_{\mathrm{conservative\,correction}}+\underbrace{\zeta\left(\boldsymbol{\mathcal{D}}(\mathbf{x})\mathbf{H}_{0}(\mathbf{x},\dot{\mathbf{x}},\tau)+\left(\mathbf{K}_{1}+\boldsymbol{\mathcal{C}}(\mathbf{x})\right)\dot{\mathbf{x}}\right)}_{\mathrm{damping\,terms}}\right]+\mathcal{O}\left(\epsilon^{2}\right)=\alpha\mathbf{q}(\tau),\label{eq:SFD-ROM-beam}
\end{gather}
where
\begin{align}
\mathbf{G}_{0}(\mathbf{x},\dot{\mathbf{x}},\tau) & =-\mathbf{K}_{2}^{-1}\boldsymbol{\mathcal{H}}(\mathbf{x})\,,\nonumber \\
\mathbf{H}_{0}(\mathbf{x},\dot{\mathbf{x}},\tau) & =-\mathbf{K}_{2}^{-1}\left[\partial_{\mathbf{x}}\boldsymbol{\mathcal{H}}(\mathbf{x})\right]\dot{\mathbf{x}}\,,\nonumber \\
\mathbf{G}_{1}(\mathbf{x},\dot{\mathbf{x}},\tau) & =-\zeta\left(\mathbf{H}_{0}(\mathbf{x},\dot{\mathbf{x}},\tau)+\mathbf{K}_{2}^{-1}\boldsymbol{\mathcal{E}}(\mathbf{x})\dot{\mathbf{x}}\right)+\mathbf{\beta K}_{2}^{-1}\mathbf{p}(\tau)\,.\label{eq:OepsTerms}
\end{align}
Note that the reduced order model (\ref{eq:SFD-ROM-beam}) is conservative
(contains only inertial and elastic force terms) at leading order,
whereas the full system (\ref{eq:FEdisc}) is dissipative due to viscoelastic
damping. The $\mathcal{O}(\epsilon)$ terms in the ROM account for
these damping contributions and hence cannot be ignored. Apart from
the damping contributions, there also exists a conservative correction
at the $\mathcal{O}(\epsilon)$ level which includes the static response
of the $\mathbf{y}$ variables to the corresponding loading $\beta\mathbf{p}(\tau)$
(cf. the expression for $\mathbf{G}_{1}(\mathbf{x},\dot{\mathbf{x}},\tau)$
in (\ref{eq:OepsTerms})).

\subsection{Specific results from SFD}

\begin{figure}[H]
\subfloat[{Displacement {[}-{]} of the beam at quarter length in the transverse
direction (slow $\mathbf{x}$ DOF).}]{\centering{}\includegraphics[width=0.48\linewidth]{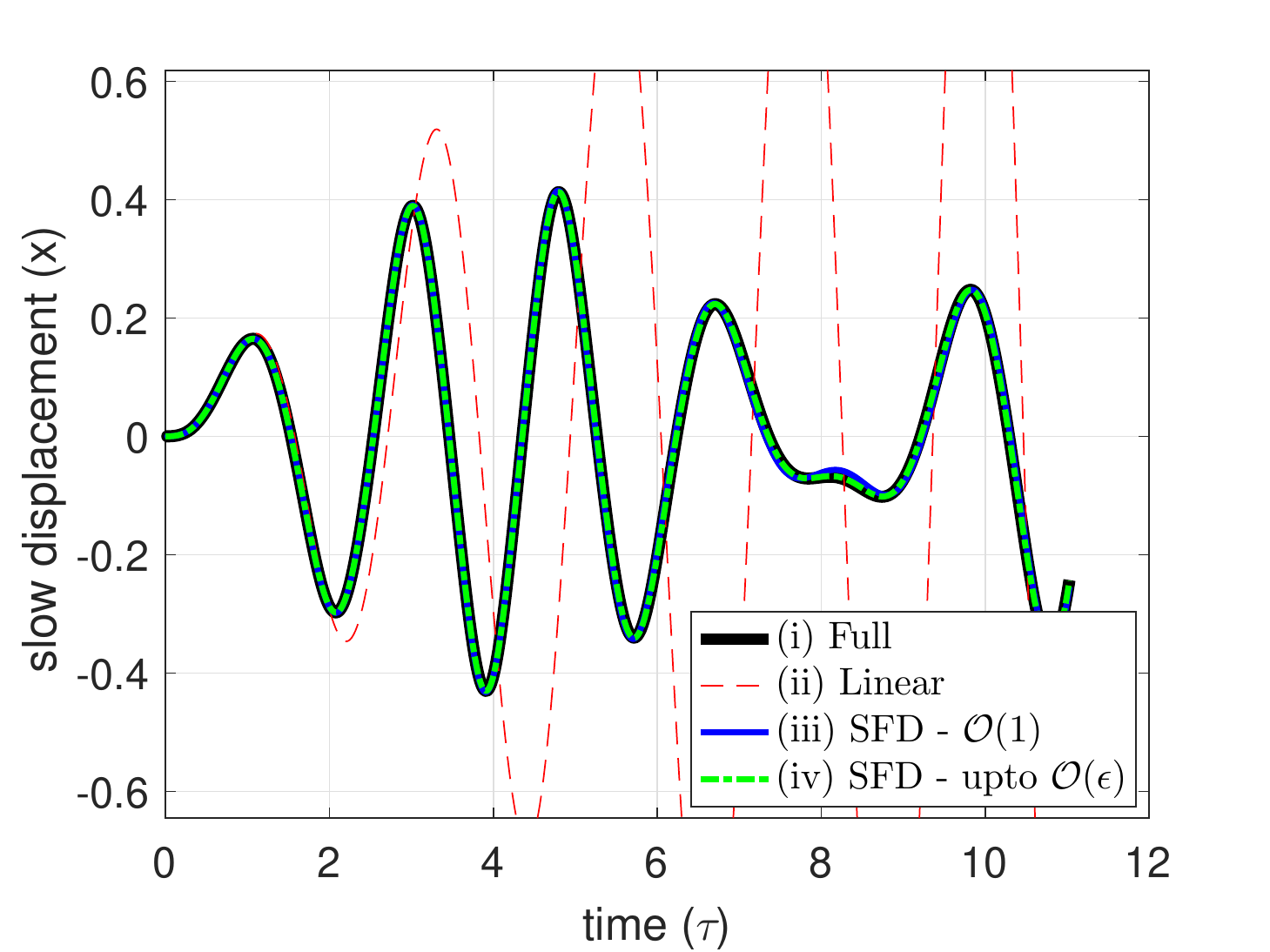}}\hfill{}\subfloat[{Displacement {[}-{]} of the beam at quarter length in the axial direction
(fast $\mathbf{y}$ DOF).}]{\centering{}\includegraphics[width=0.48\linewidth]{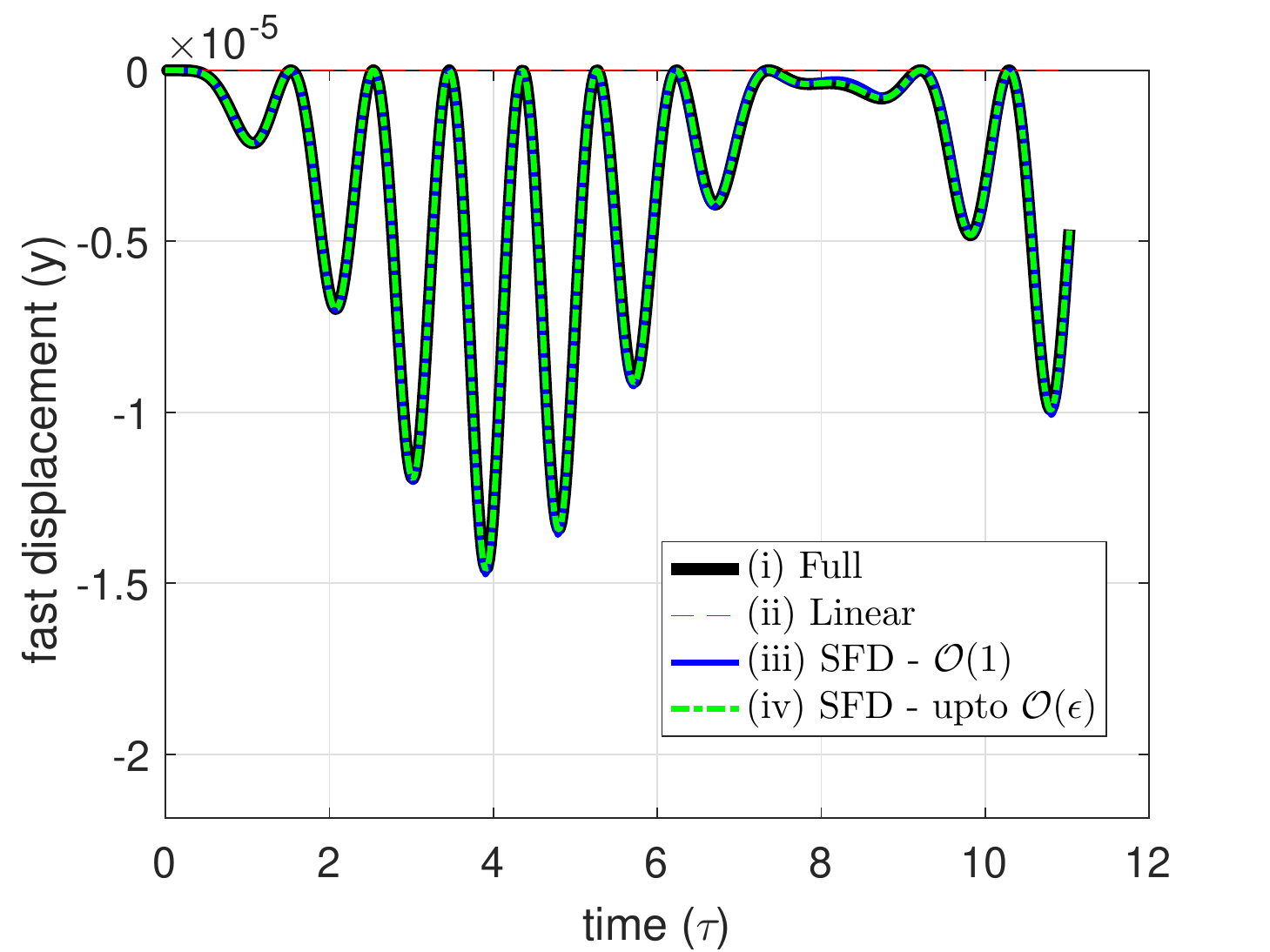}}

\caption{\label{fig:SFD-e-4}Comparison of the slow and fast components of
the reduced solution with their full nonlinear and linearized counterparts
for $\epsilon=10^{-4}$. Note that for such small values of $\epsilon$,
the ROM at the leading order (containing only conservative terms)
is a good enough representation of the full system, and is practically
identical to ROMs obtained by inclusion of the $\mathcal{O}(\epsilon)$
and $\mathcal{O}\left(\epsilon^{2}\right)$ terms. }
\end{figure}

We now consider a specific beam with geometric and material parameters
as follows: length $L=1$ m; thickness to length ratio $\epsilon$
in the range from $10^{-4}-10^{-2}$; Young's modulus $E=70$ GPa;
density$\rho=2700$ Kg/m$^{3}$; material viscous damping rate$\kappa=10^{8}$
Pa s. We use a spatially uniform load on the beam in the axial and
in the transverse direction, given by $\alpha=1$ , $\beta=1$, $q(x,\tau)=p(x,\tau)=\sin\left(\Omega T_{0}\tau\right)$.
Here $T_{0}=\frac{L}{\epsilon}\sqrt{\frac{\rho}{E}}$ is the constant
used to nondimensionalize time and $\Omega$ is the loading frequency
(chosen to be the first natural frequency of the beam in this case,
cf. Figure \ref{fig:SFD-e-4}). Using these parameters, we obtain
$\zeta\approx7.2739$.
\begin{figure}[h]
\begin{centering}
\includegraphics[width=0.6\textwidth]{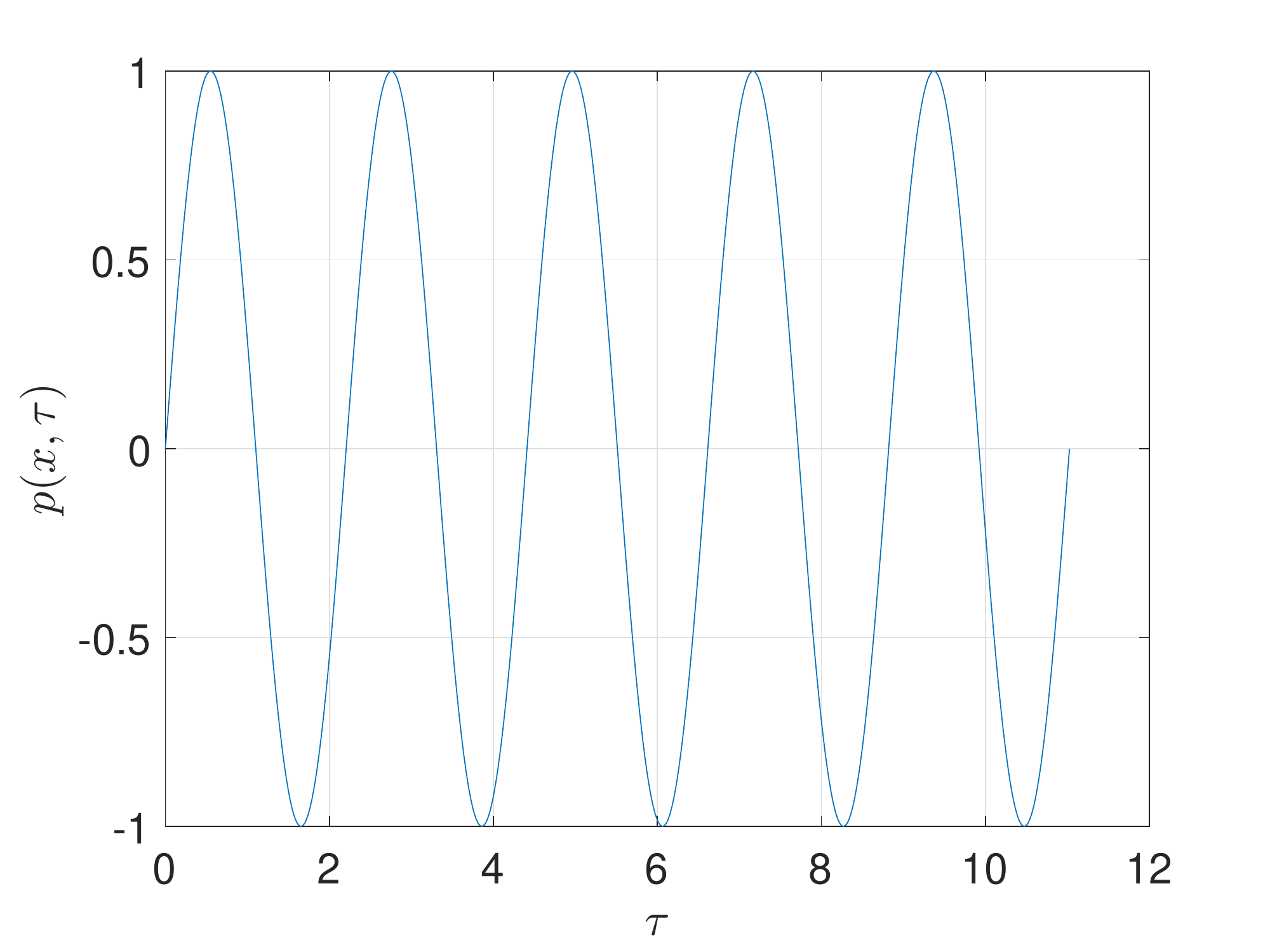}
\par\end{centering}
\caption{\label{fig:projection}The applied loading}
\end{figure}

\begin{figure}[H]
\subfloat[{Displacement {[}-{]} of the beam at quarter length in the transverse
direction (slow $\mathbf{x}$ DOF).}]{\centering{}\includegraphics[width=0.48\linewidth]{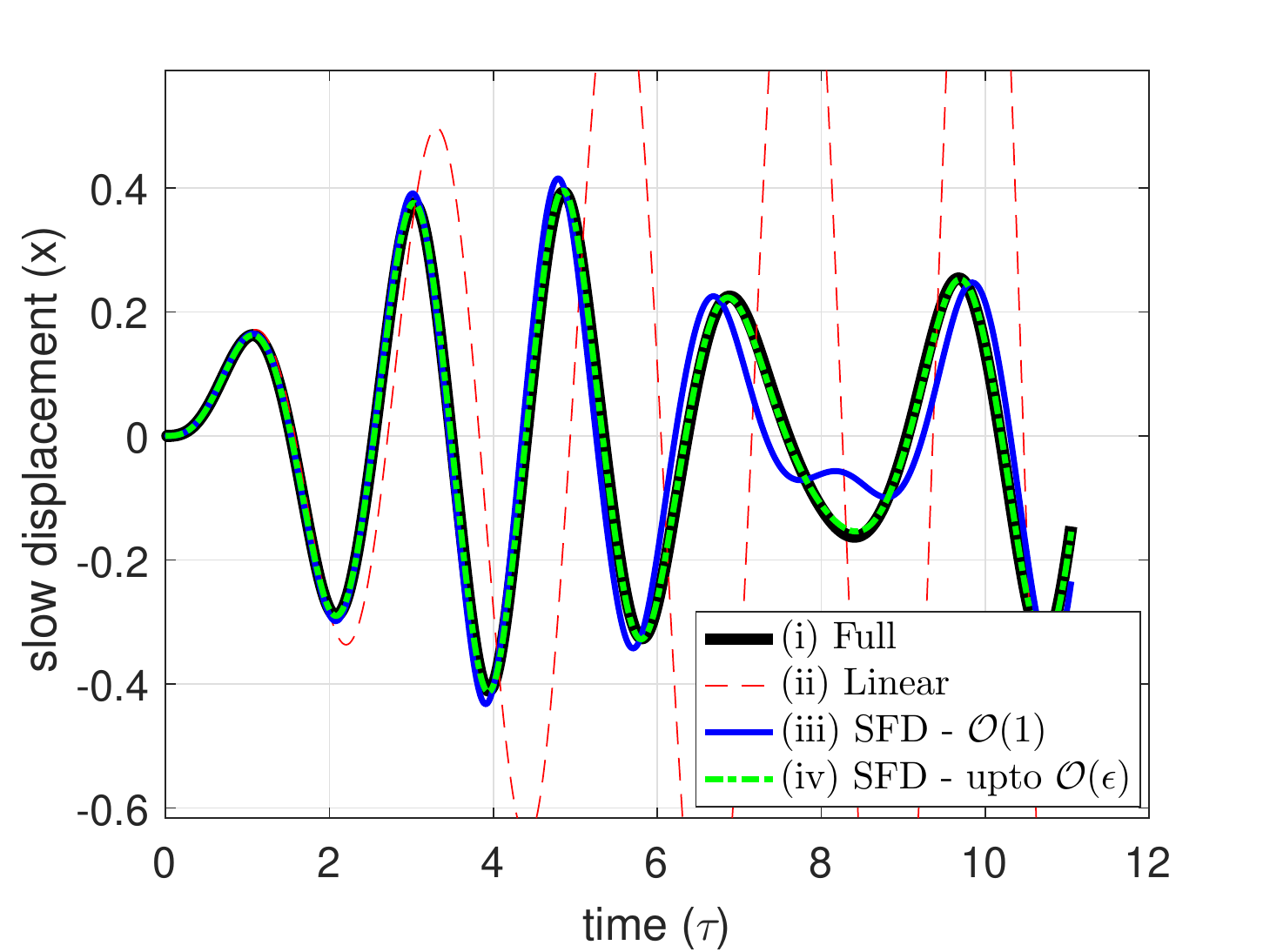}}\hfill{}\subfloat[{Displacement {[}-{]} of the beam at quarter length in the axial direction
(fast $\mathbf{y}$ DOF).}]{\centering{}\includegraphics[width=0.48\linewidth]{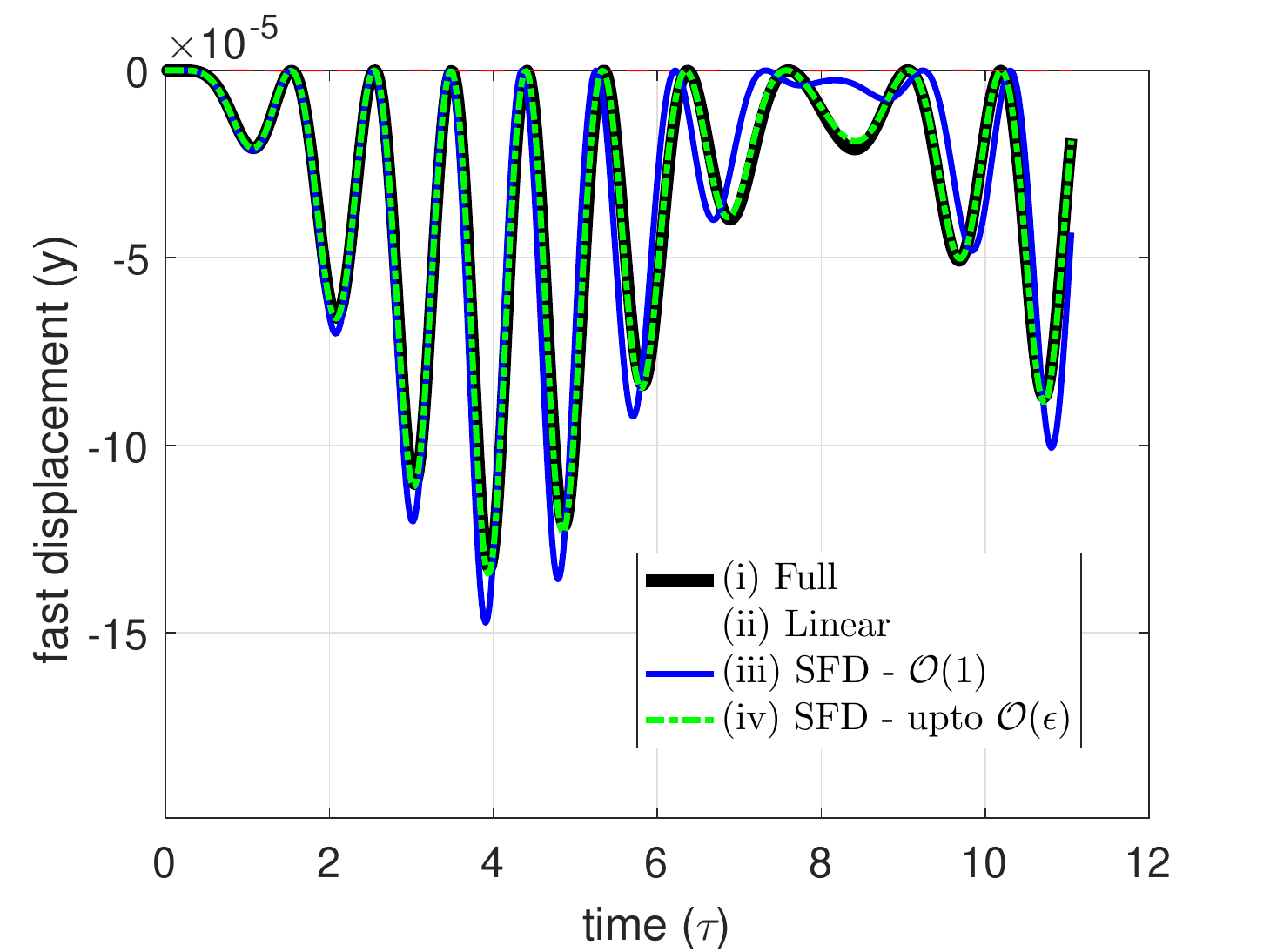}~}

\caption{\label{fig:SFD-e-3}Comparison of the slow and fast components of
the reduced solution with their full nonlinear and linearized counterparts
for $\epsilon=10^{-3}$ . For this larger value of $\epsilon$, the
leading-order conservative ROM is not accurate enough, and hence $\mathcal{O}(\epsilon)$
terms (which include damping contributions) are required to improve
accuracy. }
\end{figure}

\begin{figure}[H]
\subfloat[{Displacement {[}-{]} of the beam at quarter length in the transverse
direction (slow $\mathbf{x}$ DOF).}]{\centering{} \includegraphics[width=0.48\linewidth]{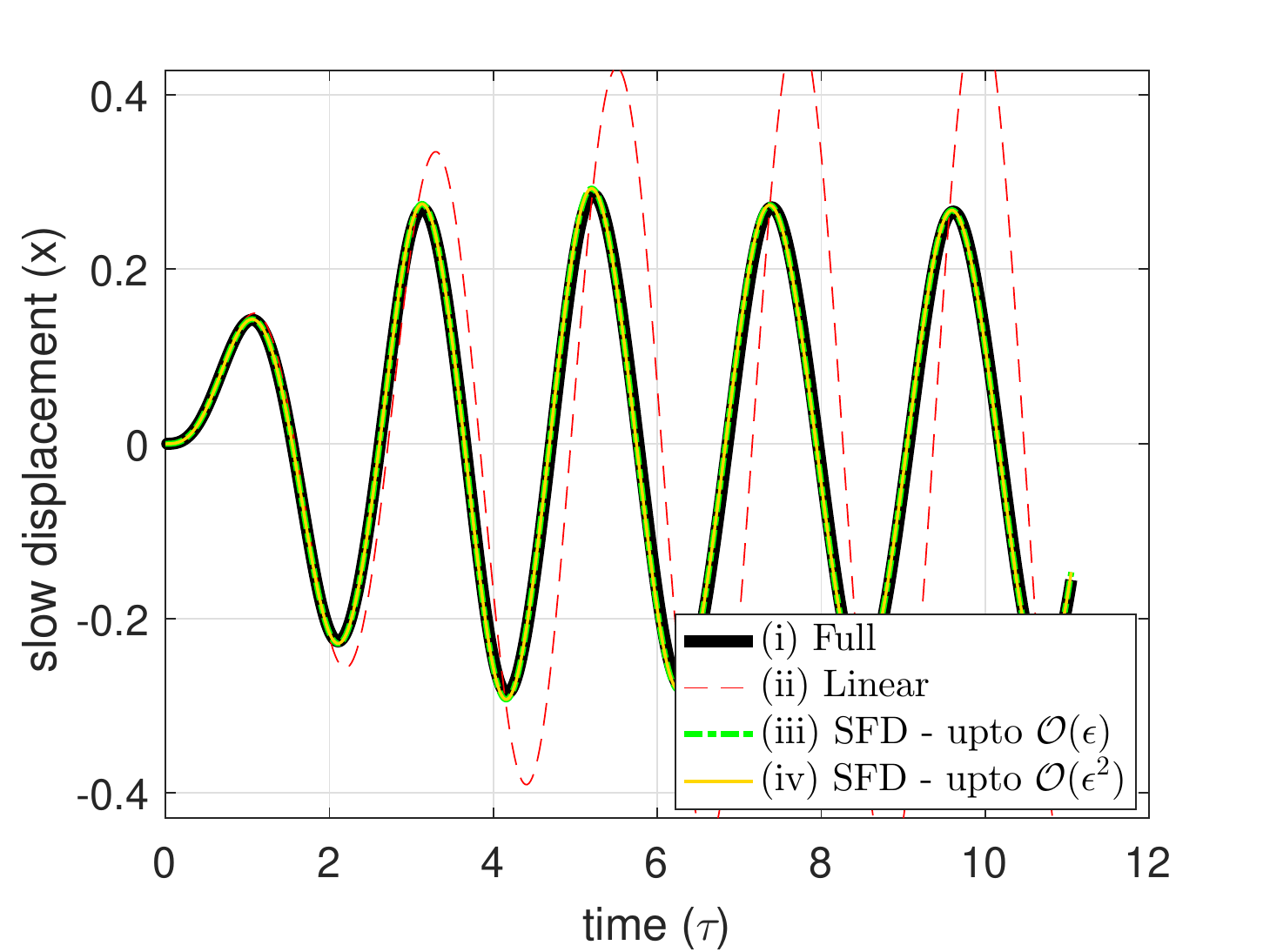}}\hfill{}\subfloat[{Displacement {[}-{]} of the beam at quarter length of beam in the
axial direction (fast $\mathbf{y}$ DOF).}]{\centering{}\includegraphics[width=0.48\linewidth]{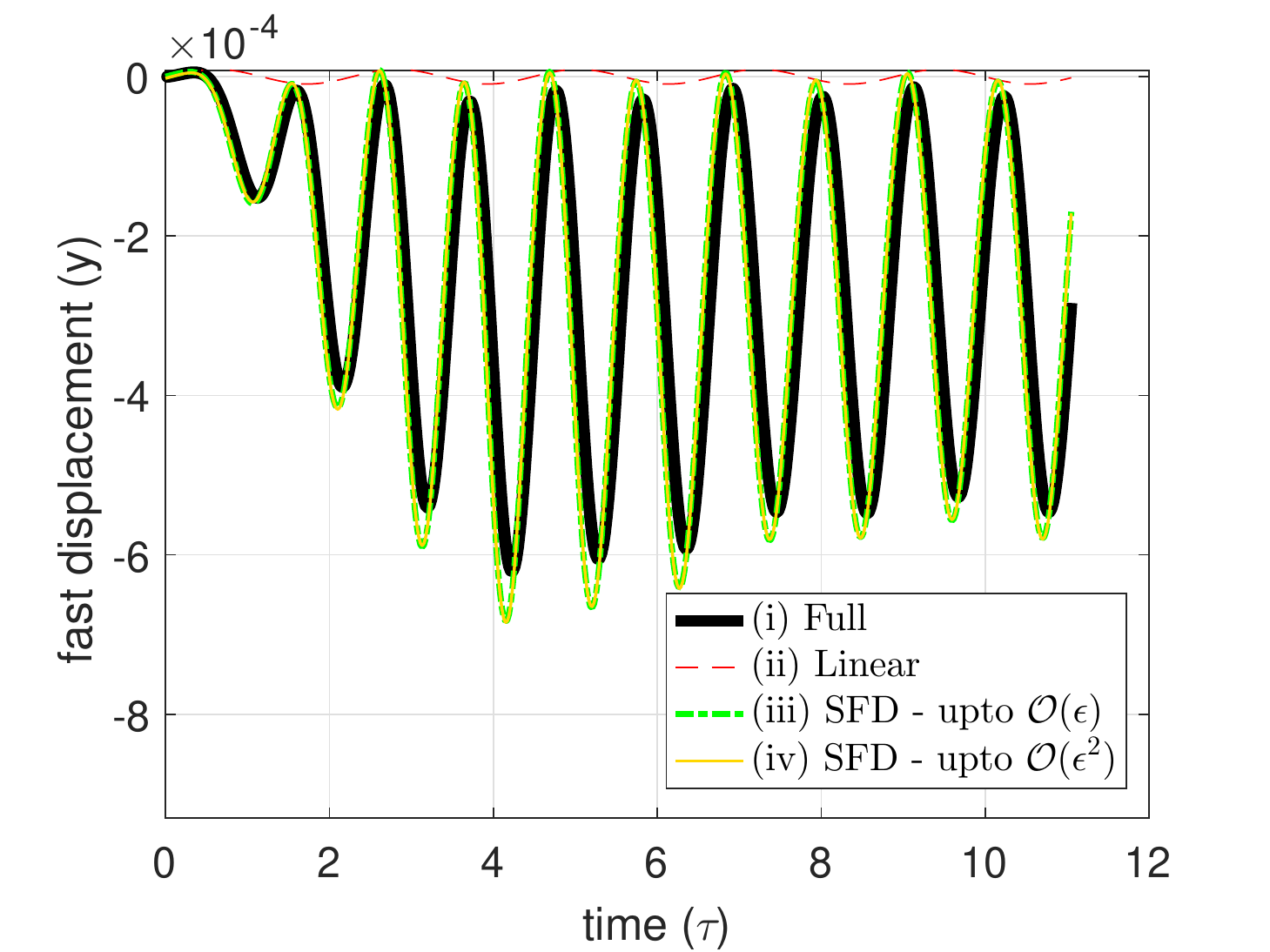}~}

\caption{\label{fig:SFD-e-2}Comparison of the slow and fast components of
the reduced solution with their full nonlinear and linearized counterparts
for $\epsilon=10^{-2}$. The $\mathcal{O}\left(\epsilon^{2}\right)$
terms were added in the ROM to improve accuracy (especially in the
slave variables), but only a marginal improvement was observed.}
\end{figure}

The graphs presented above show the response of a single degree of
freedom in time. In order to check the performance of SFD globally,
we use the following measure,
\begin{equation}
E=\dfrac{\sqrt{\sum\limits _{\tau\in S}(\mathbf{u}(\tau)-\tilde{\mathbf{u}}(\tau))^{T}(\mathbf{u}(\tau)-\tilde{\mathbf{u}}(\tau))}}{\sqrt{\sum\limits _{\tau\in S}\mathbf{u}(\tau)^{T}\mathbf{u}(\tau)}}\times100\%\,,\label{eq:errorMeasure}
\end{equation}
where $\mathbf{u}(\tau)\in\mathbb{R}^{n_{s}+n_{f}}$ is the full vector
of generalized displacements at time $\tau$, obtained from the full
nonlinear solution; $\tilde{\mathbf{u}}(\tau)\in\mathbb{R}^{n_{s}+n_{f}}$
is the solution based on the reduced model; and $S$ is the set of
time instants at which the error is recorded. The error recorded for
different cases is shown in Table \ref{tab:Relative-Error}. We observe
that for $\epsilon$ of order $10^{-4}$, the ROM at the leading order
(which is conservative) provides a close approximation for the full
system. The $\mathcal{O}(\epsilon)$ terms, however, become important
when $\epsilon$ is of order $10^{-3}$, because they account for
damping which becomes more significant in this$\epsilon$-range. We
further calculated the $\mathcal{O}\left(\epsilon^{2}\right)$ terms
in the ROM to compare changes in accuracy (see \ref{sec:appA} for
general calculation of the $\mathcal{O}\left(\epsilon^{2}\right)$
terms). We observe that the $\mathcal{O}\left(\epsilon^{2}\right)$
terms only provide a marginal improvement in terms of accuracy in
our case. A further increase in $\epsilon$ values to the order of
$10^{-1}$ renders the reduced model highly inaccurate, thus higher-order
terms in $\epsilon$ are necessary to achieve the desired accuracy.
Note, however, the $\epsilon$ values around $10^{-1}$ are rather
unphysical; the parameter $\epsilon$, denoting the thickness-to-length
ratio of a \textit{beam}, is not expected attain such high values. 

\begin{table}[H]
\begin{centering}
\begin{tabular}{|c|c|c|c|}
\hline 
$\epsilon$ & SFD - $\mathcal{O}(1)$ & SFD - up to $\mathcal{O}(\epsilon)$ & SFD - up to $\mathcal{O}\left(\epsilon^{2}\right)$\tabularnewline
\hline 
\hline 
$10^{-4}$ & 3.7078 & 0.3670 & 0.3670\tabularnewline
\hline 
$10^{-3}$ & 29.589 & 2.9681 & 2.9681\tabularnewline
\hline 
$10^{-2}$ & >100 & 5.3432 & 5.3425\tabularnewline
\hline 
$10^{-1}$ & >100 & >100 & >100\tabularnewline
\hline 
\end{tabular}\caption{\label{tab:Relative-Error}Relative reduction error $E$ calculated
according from formula (\ref{eq:errorMeasure}) }
\par\end{centering}
\end{table}

\section{SSM reduction}

Though the SFD reduction is robust and globally valid in the phase
space, the reduced system on the slow manifold still contains $n_{s}$
degrees of freedom, which is generally a large number. Specifically,
due to the choice of the shape functions for discretization, the reduced
model (\ref{eq:SFD-ROM-beam}) obtained from SFD still contains two-thirds
of the total number of unknowns of the full system (\ref{eq:FEdisc}).
The following eigenvalue analysis of the linearized reduced system,
however, reveals the existence of a further separation in time scales
within the slow manifold near the origin, showing potential for a
further reduction via spectral submanifolds (cf. Haller and Ponsioen
\cite{haller16})

\subsection{Separation in damping: spectral quotient analysis}

The eigenvalue problem for the linearization of the reduced system
(\ref{eq:SFD-ROM-beam}) can be formulated as 
\begin{equation}
(\lambda_{k}^{2}\mathbf{M}_{1}+\lambda_{k}\zeta\epsilon\mathbf{K}_{1}+\mathbf{K}_{1})\boldsymbol{\phi}_{k}=\mathbf{0}\,,\label{eq:EVP}
\end{equation}
where $\lambda_{k}\in\mathbb{C}$ and $\boldsymbol{\phi}_{k}\in\mathbb{C}^{n_{s}}$
are the $k^{\mathrm{th}}$ eigenvalue and eigenvector for any $k\in\{1,\dots,n_{s}\}$.
The negative real part of the eigenvalue $\lambda_{k}$ represents
the exponential rate of decay of trajectories towards the equilibrium
position along the corresponding two-dimensional eigenspace of the
linearized system. The smoothest local extension of such an invariant
subspace is known as the spectral submanifold (SSM), a notion introduced
by Haller and Ponsioen \cite{haller16}. Tangent to a slow subspace
(i.e.$\,$, the subspace spanned by eigenvectors corresponding to
eigenvalues with the lowest-magnitude real parts), such an SSM offers
an opportunity for a further, drastic reduction from the SFD-based
slow manifold to a two-dimensional invariant manifold. 

\begin{figure}[H]
\begin{centering}
\subfloat[]{\centering{} \includegraphics[width=0.45\linewidth]{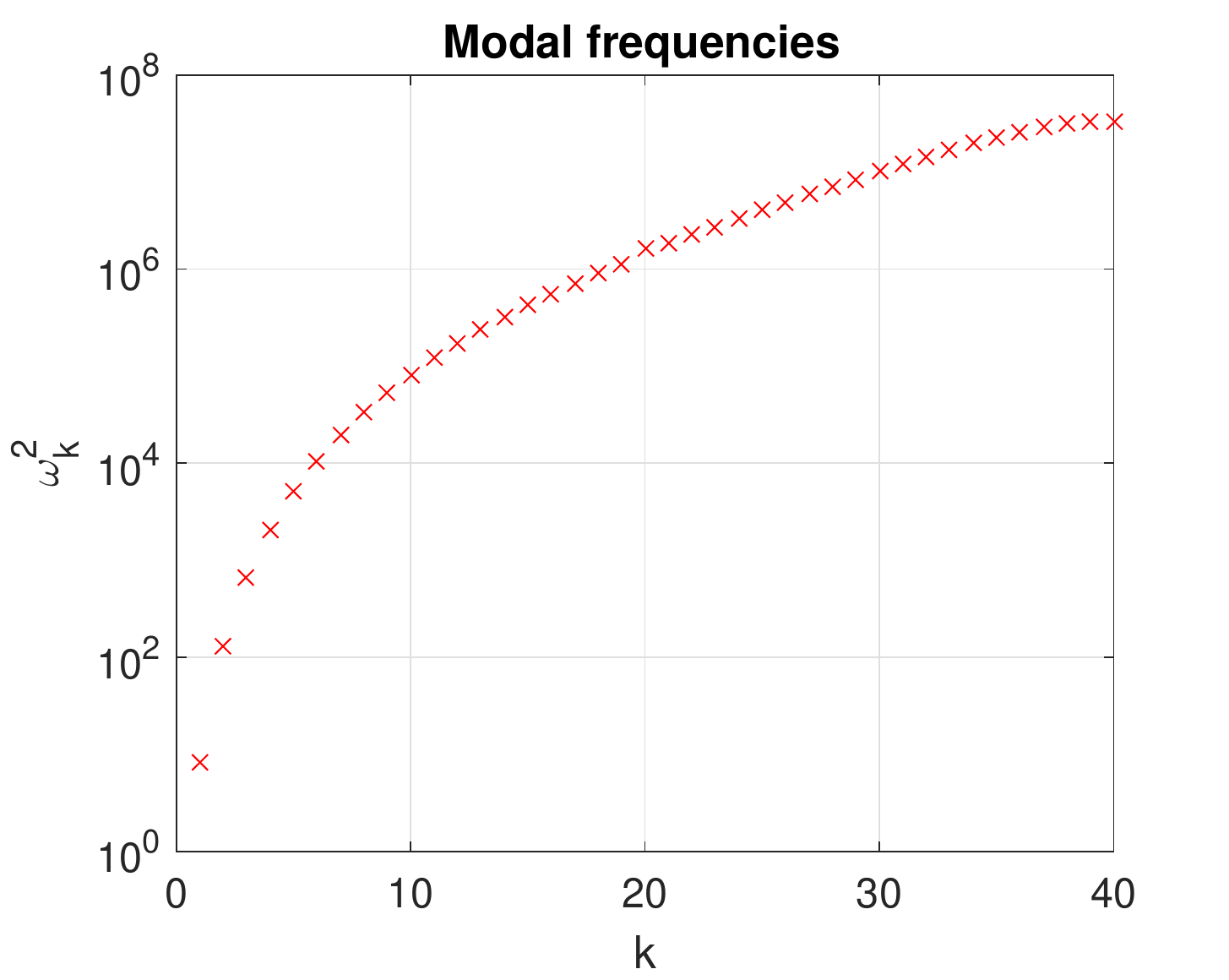}}\subfloat[]{\centering{}\includegraphics[width=0.45\linewidth]{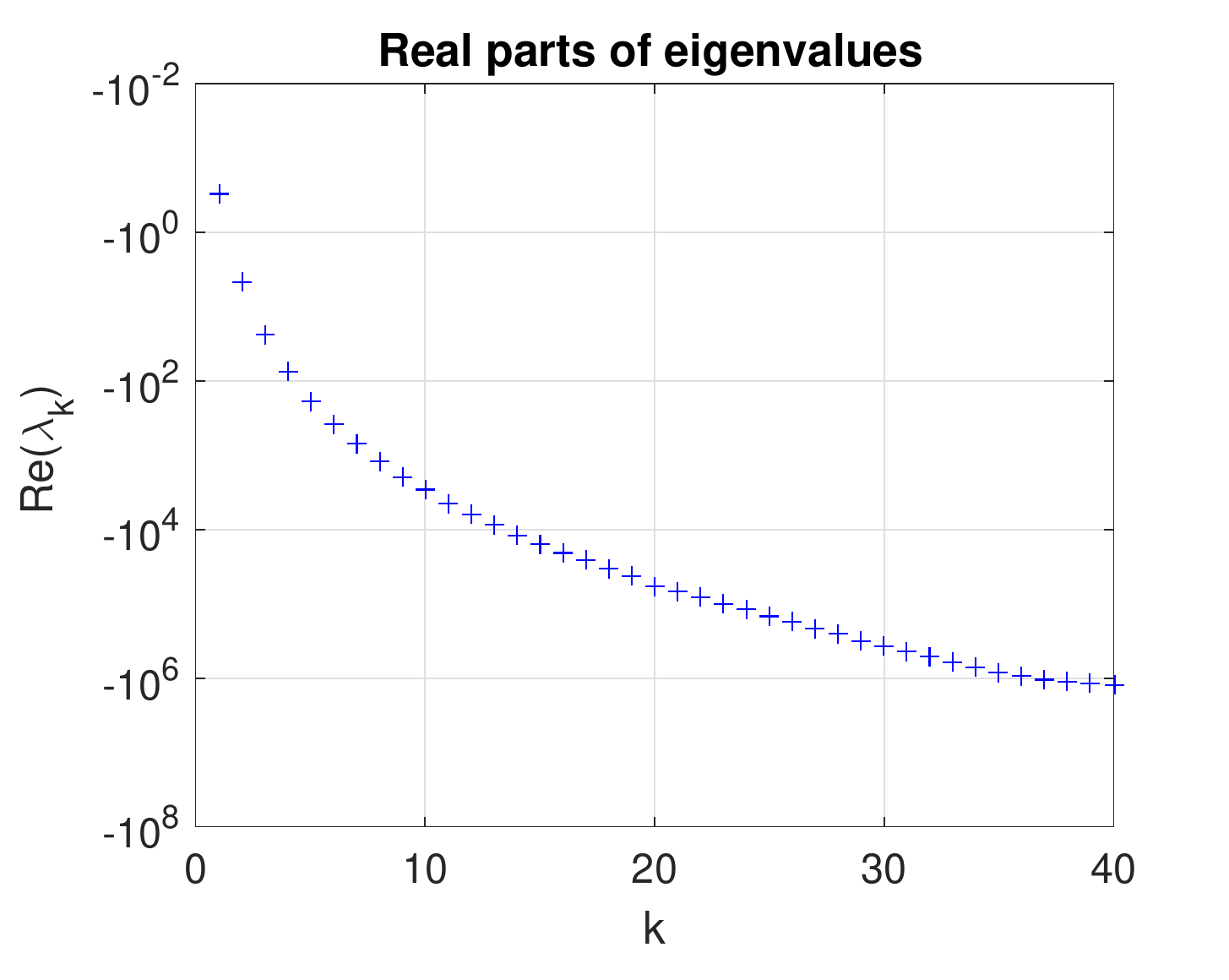}~}
\par\end{centering}
\centering{}\subfloat[]{\centering{}\includegraphics[width=0.45\linewidth]{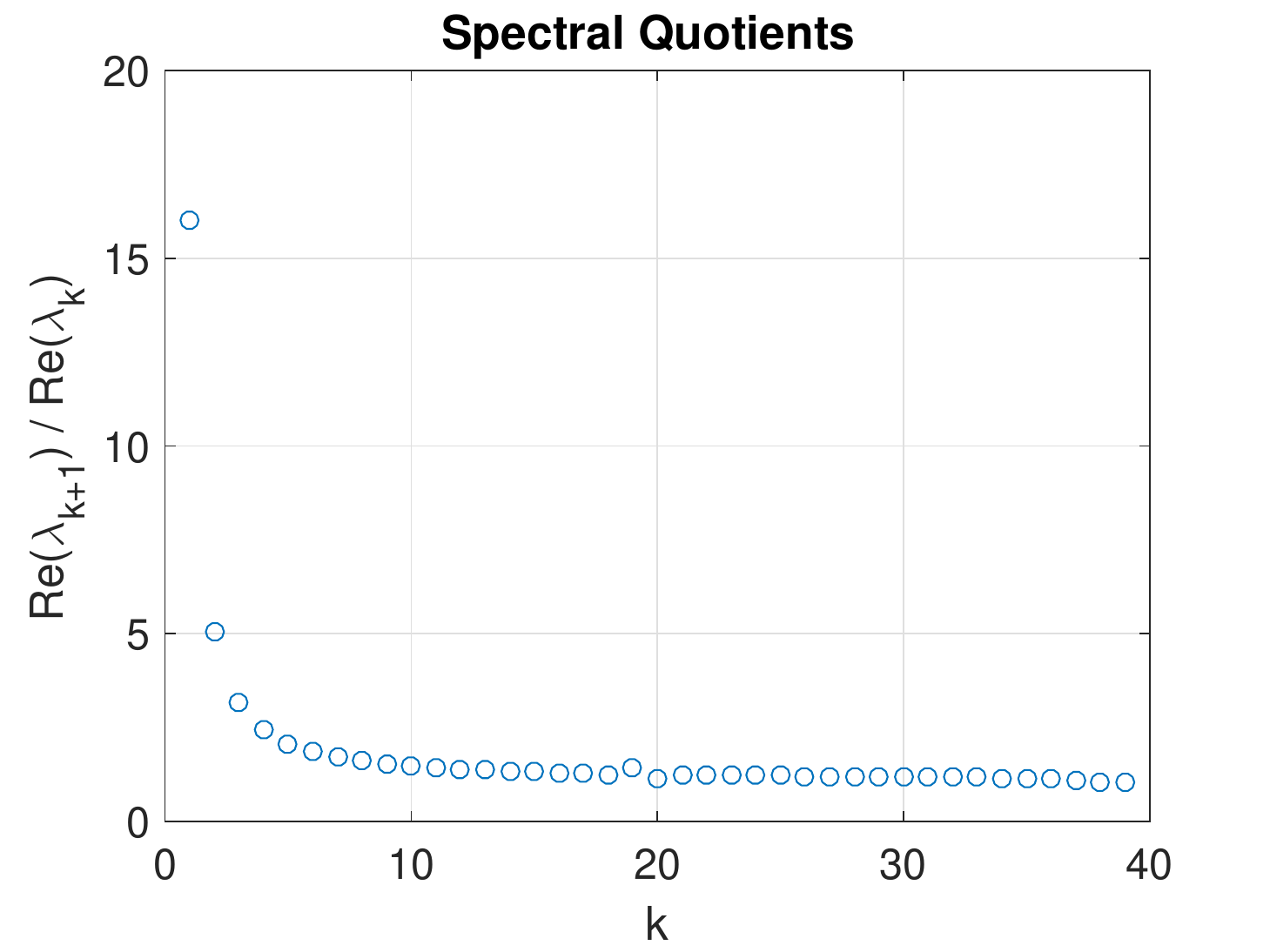}~}\caption{\label{fig:DampingGaps}Spectral quotient analysis: the approximate
ratio (\ref{eq:EigenvalueApproximation}) between the real parts of
successive eigenvalues of the system (c) reveals the spectral quotient
gaps in the system. The first few eigen-modes form a subspace much
slower than the rest and hence are relavant for the SSM reduction.
No such large spectral gaps can be identified from the plots of the
undamped natural frequencies (a) and real parts of the eigenvalues
(b). }
\end{figure}

Estimating the real parts of the eigenvalues arising from (\ref{eq:EVP})
is sensitive to the choice of the eigensolver algorithm, especially
for large $n_{s}$ values. An accurate first-order approximation to
the eigenvalues, however, can be obtained using the undamped eigenvalues
from their defining equation 
\[
(\mathbf{K}_{1}-\omega_{0k}^{2}\mathbf{M}_{1})\boldsymbol{\phi}_{0k}=\mathbf{0}\,.
\]
 Specifically, for small damping and well-separated eigenvalues, a
reliable first-order approximation for the real parts of the eigenvalues
of system (\ref{eq:EVP}) can be obtained as follows \cite{geradin15}:
\begin{equation}
\lambda_{k}=-\frac{\boldsymbol{\phi}_{0k}^{T}(\zeta\epsilon\mathbf{K}_{1})\boldsymbol{\phi}_{0k}}{2\boldsymbol{\phi}_{0k}^{T}\mathbf{M}_{1}\boldsymbol{\phi}_{0k}}+i\,\omega_{0k}^{2}+\mathcal{O}\left(\epsilon^{2}\right)=-\frac{\zeta\epsilon}{2}\omega_{0k}^{2}+i\,\omega_{0k}^{2}+\mathcal{O}\left(\epsilon^{2}\right)\quad\forall k\in\{1,\dots,n_{s}\}\,.\label{eq:EigenvalueApproximation}
\end{equation}
Based on this formula, the modal frequencies, real parts of the eigenvalues,
and the ratios of the subsequent real parts are shown as a function
of the mode number in Figure \ref{fig:DampingGaps}.
 Note that the first eigenspace is about $16$ times slower in terms
of its exponential decay rate than the second slowest eigenspace,
which renders the first mode an optimal choice for reduction using
its corresponding SSM. We will, therefore, use this mode to perform
a single-mode SSM reduction for the beam model, as discussed below. 

\subsection{SSM-based reduction }

In order to perform a further reduction of the reduced system (\ref{eq:SFD-ROM-beam})
obtained from SFD, we now compute the slowest single-mode SSM, whose
reduced dynamics is given by a two-dimensional ODE as a final reduced-order
model. Szalai et al. \cite{Szalai16} have obtained general formulae
for such one-mode SSMs and their associated backbone curves for autonoumous
dynamical systems (no external forcing in our setting). The reduced-order
model obtained in this fashion is given in terms of an internal parametrization
of the SSM, rather than a projection of the SSM on its underlying
modal subspace. This construct allows the SSM to be recovered more
globally, even if it develops a fold over the underlying modal subspace. 

For this general formulation, we require the autonomous counterpart
of system (\ref{eq:SFD-ROM-beam}) to be in the phase space form 
\begin{align}
\underbrace{\left[\begin{array}{c}
\dot{\mathbf{x}}\\
\ddot{\mathbf{x}}
\end{array}\right]}_{\dot{\mathbf{\tilde{z}}}} & =\underbrace{\left[\begin{array}{cc}
\mathbf{0} & \mathbf{I}\\
-\mathbf{M}_{1}^{-1}\mathbf{K}_{1} & -\epsilon\zeta\mathbf{M}_{1}^{-1}\mathbf{K}_{1}
\end{array}\right]}_{\mathbf{A}}\underbrace{\left[\begin{array}{c}
\mathbf{x}\\
\dot{\mathbf{x}}
\end{array}\right]}_{\mathbf{\tilde{z}}}+\mathbf{t}(\mathbf{x},\dot{\mathbf{x}})\,,\nonumber \\
\mathrm{or}\quad\dot{\mathbf{\tilde{z}}} & =\mathbf{A}\tilde{\mathbf{z}}+\mathbf{t}(\tilde{\mathbf{z}})\,,\label{eq:statespaceROM}
\end{align}
where $\mathbf{t}(\mathbf{x},\dot{\mathbf{x}})$ is a class $C^{r}$
function representing the nonlinear terms in the ROM (\ref{eq:SFD-ROM-beam}).
By inspection of (\ref{eq:SFD-ROM-beam}), we find that $\mathbf{t}$
is a strictly cubic polynomial in $\mathbf{\tilde{z}}$ and hence
belongs to the class $C^{a}$ of analytic functions. The system (\ref{eq:statespaceROM})
can be diagonalized as 
\begin{equation}
\dot{\mathbf{z}}=\mathbf{\boldsymbol{\Lambda}}\mathbf{z}+\boldsymbol{\mathcal{T}}\left(\mathbf{z}\right),\quad\mathbf{z}\in\mathbb{C}^{2n_{s}},\quad\mathbf{\Lambda}=\mathrm{diag}\left(\lambda_{1},\dots,\,\lambda_{2n_{s}}\right),\,\lambda_{2i-1}=\bar{\lambda}_{2i},\,i\in\{1,\dots,\,n_{s}\},\quad\boldsymbol{\mathcal{T}}(\mathbf{z})=\mathcal{O}\left(|\mathbf{z}|^{3}\right)\,,\label{eq:diagROM}
\end{equation}
where $\lambda_{i}\in\mathbb{C\,}$ are the eigenvalues of the matrix
$\mathbf{A}$ ordered such that $\mathrm{Re}\left(\lambda_{n_{s}}\right)\le\dots\le\mathrm{Re}\left(\lambda_{1}\right)<0$.
The matrix $\mathbf{P}$ of the linear transformation $\mathbf{z}=\mathbf{P}^{-1}\tilde{\mathbf{z}}$
leading to (\ref{eq:diagROM}), contains the corresponding eigenvectors
of $\mathbf{A}$. 

The notion of an SSM was introduced by Haller and Ponsioen \cite{haller16},
more formally defined as follows:
\begin{defn}
A \textit{spectral submanifold} (SSM), $W(\mathcal{E})$, corresponding
to a spectral subspace $\mathcal{E}$ of $\mathbf{\Lambda}$ is an
invariant manifold of the dynamical system (\ref{eq:diagROM}) such
that 
\end{defn}
\begin{elabeling}{00.00.0000}
\item [{(i)}] $W(\mathcal{E})$ is tangent to $\mathcal{E}$ at the origin
and has the same dimension as $\mathcal{E}$;
\item [{(ii)}] $W(\mathcal{E})$ is strictly smoother than any other invariant
manifold satisfying (i).
\end{elabeling}
When $\boldsymbol{\mathcal{T}}\left(\mathbf{z}\right)$ is a $C^{r}$
function, then under the assumptions 
\begin{description}
\item [{(B1)}] the \textit{relative spectral quotient} $\sigma(\mathcal{E}):=\mathrm{Int}\left[\frac{\min_{j\ne\ell,\ell+1}\mathrm{Re}(\lambda_{j})}{\mathrm{Re}(\lambda_{\ell})}\right]$
satisfies $\sigma(\mathcal{E})\le r$, 
\item [{(B2)}] there are no resonances up to order $\sigma(\mathcal{E})$
between $\lambda_{\ell},\,\bar{\lambda}_{\ell}$ and the rest of the
spectrum of $\mathbf{\Lambda}$,
\end{description}
the existence and uniqueness of a class $C^{r+1}$ SSM is guaranteed
by the main theorem of Haller and Ponsioen \cite{haller16}, deduced
from the more general results of Cabré et al. \cite{cabre2003}. The
theorem states that the SSM can be viewed as an embedding of an open
set $\mathcal{U}\subset\mathbb{C}^{2}$ into the phase space of system
(\ref{eq:diagROM}), described by a map $\mathbf{W}:\mathcal{U}\subset\mathbb{C}^{2}\mapsto\mathbb{C}^{2n_{s}}$.
Furthermore, there exists a quadratic polynomial mapping $\mathbf{R}:\mathcal{U}\mapsto\mathcal{U}$,
such that the reduced dynamics on the SSM can be written as 
\begin{equation}
\dot{\mathbf{s}}=\mathbf{R\left(s\right)},\qquad\mathbf{R\left(s\right)}=\begin{bmatrix}\lambda_{\ell}z_{\ell}+\beta_{\ell}z_{\ell}^{2}\bar{z}_{\ell}\\
\bar{\lambda}_{\ell}\bar{z}_{\ell}+\bar{\beta}_{\ell}z_{\ell}\bar{z}_{\ell}^{2}
\end{bmatrix},\label{eq:SSMROM}
\end{equation}
where $\mathbf{s}=(z_{\ell},\bar{z}_{\ell})$ are the coordinates
for the mode $\ell$ (for which the SSM construction is performed). 

In our beam setting, we obtain $\sigma(\mathcal{E})=\frac{\omega_{0n_{s}}^{2}}{\omega_{0\ell}^{2}}+\mathcal{O}(\epsilon)$
from the eigenvalue approximation in (\ref{eq:EigenvalueApproximation}).
This shows that the relative spectral quotient $\sigma(\mathcal{E})$
depends on the chosen discretization and monotonically increases with
the number of discretization variables. This can be expected physically
from a proportionally damped mechanical structure in which the damping
for the modes is expected to increase monotonically with the corresponding
oscillation frequency. Since $\mathbf{g}\in C^{a}$ here, we obtain
that (\textbf{B1}) is satisfied. Furthermore, we assume generic parameter
values under which the non-resonance requirement (\textbf{B2}) is
satisfied by the beam system. The theorem of Haller and Ponsioen \cite{haller16}
then applies to the system (\ref{eq:diagROM}), and the formulae obtained
by Szalai et al. \cite{Szalai16} yield the coefficients of $\mathbf{R}$
and the quadratic and cubic Taylor coefficients of $\mathbf{W}$.
Specifically, let $\boldsymbol{\mathcal{T}}$ be expanded in indicial
notation as 
\begin{equation}
\left(\mathbf{\boldsymbol{\mathcal{T}}}(\mathbf{z})\right)_{i}=\sum_{j=1}^{2n_{s}}\sum_{k=1}^{2n_{s}}\sum_{l=1}^{2n_{s}}T_{ijkl}z_{j}z_{k}z_{l}\quad\forall i,j,k,l\in{1,\dots,2n_{s}}\,,\label{eq:Texpansion}
\end{equation}
and let $\mathbf{W}(\mathbf{s})$ be of the polynomial form 
\begin{equation}
\mathbf{W}(\mathbf{s})=\mathbf{\boldsymbol{\mathcal{W}}}^{(1)}(\mathbf{s})+\mathbf{\boldsymbol{\mathcal{W}}}^{(2)}(\mathbf{s})+\boldsymbol{\mathcal{W}}^{(3)}(\mathbf{s})+\dots,\label{eq:Wexpansion}
\end{equation}
where $\mathbf{\boldsymbol{\mathcal{W}}}^{(n)}(\mathbf{s})$ denotes
the $n^{\mathrm{th}}$order terms of $\mathbf{W}$. 

The first-order terms in this expansion are given by 
\[
\mathbf{\boldsymbol{\mathcal{W}}}^{(1)}(\mathbf{s})=\mathbf{W}^{(1)}\mathbf{s}\,,
\]
where $\mathbf{W}^{(1)}\in\mathbb{C}^{2n_{s}\times2}$ is an all-zero
matrix except for two non-zero entries given by $\left(\mathbf{W}^{(1)}\right)_{1,\ell}=\lambda_{\ell},\,\left(\mathbf{W}^{(1)}\right)_{2,\ell+1}=\lambda_{\ell+1}$.
The $i^{\mathrm{th}}$ component of the quadratic terms can be written
as
\[
\left(\mathbf{\boldsymbol{\mathcal{W}}}^{(2)}(\mathbf{s})\right)_{i}=\sum_{j=1}^{2}\sum_{k=1}^{2}W_{ijk}^{(2)}s_{j}s_{k}\equiv0\,,\quad i\in\{1,\dots,2n_{s}\},\quad j,k\in\{1,2\}\,,
\]
where $\mathbf{W}^{(2)}\in\mathbb{C}^{2n_{s}\times2\times2}$ is an
all-zero 3-tensor since $\boldsymbol{\mathcal{T}}$ has no quadratic
components. Finally, using the formulae derived by Szalai et al. \cite{Szalai16},
we can write the cubic terms as 
\[
\left(\mathbf{\boldsymbol{\mathcal{W}}}^{(3)}(\mathbf{s})\right)_{i}=\sum_{j=1}^{2}\sum_{k=1}^{2}\sum_{l=1}^{2}W_{ijkl}^{(3)}s_{j}s_{k}s_{l}\,,\quad i\in\{1,\dots,2n_{s}\},\quad j,k,l\in\{1,2\}\,,
\]
where $\mathbf{W}^{(3)}\in\mathbb{C}^{2n_{s}\times2\times2\times2}$
is a sparse 4-tensor with nonzero entries given as
\begin{align*}
W_{i111}^{(3)} & =\frac{T_{i\ell\ell\ell}}{3\lambda_{\ell}-\lambda_{i}}\,,\quad i\in\{1,\dots,2n_{s}\}\,,\\
W_{i222}^{(3)} & =\frac{T_{i(\ell+1)(\ell+1)(\ell+1)}}{3\lambda_{\ell+1}-\lambda_{i}}\,,\quad i\in\{1,\dots,2n_{s}\}\,,\\
W_{ijkl}^{(3)} & =(1-\delta_{i\ell})\frac{T_{ijkl}}{2\lambda_{\ell}+\bar{\lambda}_{\ell}-\lambda_{i}}\,,\quad i\in\{1,\dots,2n_{s}\},\quad(j,k,l)\in\{(\ell,\ell+1,\ell),(\ell,\ell,\ell+1),(\ell+1,\ell,\ell)\}\,,\\
W_{ijkl}^{(3)} & =(1-\delta_{i(\ell+1)})\frac{T_{ijkl}}{2\lambda_{\ell+1}+\bar{\lambda}_{\ell+1}-\lambda_{i}}\,\\
 & \quad i\in\{1,\dots,2n_{s}\},\quad(j,k,l)\in\{(\ell+1,\ell,\ell+1),(\ell+1,\ell+1,\ell),(\ell,\ell+1,\ell+1)\}\,.
\end{align*}
Here $T_{ijkl}\in\mathbb{C}$ denotes the components of the 4-tensor
in the expansion of $\boldsymbol{\mathcal{T}}$ given in (\ref{eq:Texpansion}),
and $\delta_{ij}$ respresents the Kronecker-delta. The expression
for $\beta_{\ell}$ used in the expansion of $\mathbf{R}$ as in (\ref{eq:SSMROM})
is given by
\[
\beta_{\ell}=T_{\ell\ell\ell(\ell+1)}+T_{\ell\ell(\ell+1)\ell}+T_{\ell(\ell+1)\ell\ell},
\]
as obtained from the general formulae of Szalai et al. \cite{Szalai16},
reproduced in the \ref{sec:appB}. Note that the Einstein summation
convention has \textit{not} been followed in the above expressions.
We finally express the reduced dynamics in the polar coordinates $\left(\rho,\theta\right)$
using a transformation $\boldsymbol{\mathcal{R}}:\mathbb{\mathbb{C}}^{2}\mapsto\mathbb{R}\times\mathbb{S}^{1}$
such that $\mathbf{s}=\boldsymbol{\mathcal{R}}^{-1}\left(\rho,\theta\right)=\left(\rho e^{i\theta},\rho e^{-i\theta}\right)$,
given by 
\begin{align}
\dot{\rho} & =\rho(\mathrm{Re}\lambda_{\ell}+\mathrm{Re}\beta_{\ell}\rho^{2})\,,\nonumber \\
\dot{\theta} & =\mathrm{Im}\lambda_{\ell}+\mathrm{Im}\beta_{\ell}\rho^{2}\,.\label{eq:polarROM}
\end{align}

\subsection{Results}

\begin{figure}[h]
\begin{centering}
\includegraphics[width=0.8\textwidth]{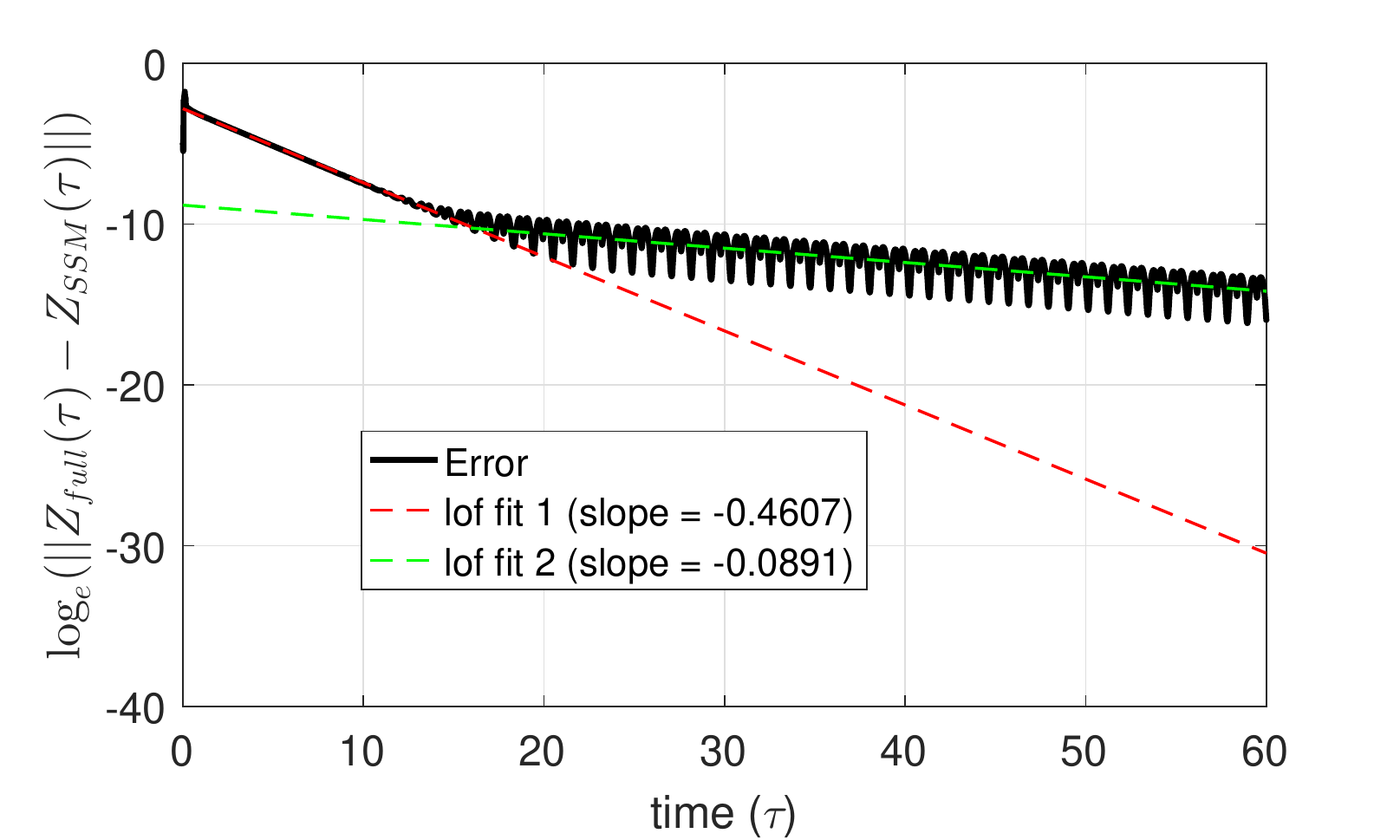}\caption{\label{fig:decayRate}The logarithmic plot shows $|\mathbf{z}(\tau)-\mathbf{W}(\mathbf{s}(\tau))|$
vs. time. The two fits suggest initial and final decay rates of $e^{-0.4607\tau}$
and $e^{-0.0891\tau}$ approximately, respectively. The approximate
decay rates along and transverse to the SSM near the equilibrium are
given by $e^{-0.0295\tau}$ and $e^{-0.4724\tau}$, respectively (based
on the eigenvalues of the linearized system).}
\par\end{centering}
\end{figure}

We now consider the SFD-reduced beam system (\ref{eq:SFD-ROM-beam})
with the physically relevant parameter value $\epsilon=10^{-3}$,
and further reduce it to its slowest, two-dimensional SSM with $\ell=1$.
Let the modal coordinates be partitioned such that $Q_{1}(\mathbf{z}(\tau))=(z_{1}(\tau),\,z_{2}(\tau))$
represents the displacement of master modes, and $Q_{2}(\mathbf{z}(\tau))=(z_{3},\dots,z_{2n_{s}})$
represents the rest of the slow modes obtained previously from the
SFD. For the two-dimensional SSM over longer time-scales, free oscillations
of the beam are expected to decay towards the origin. To illustrate
this over longer time scales, we consider three sets of intial conditions. 
\begin{enumerate}
\item \textbf{Initial condition on the SSM}: Starting with an initial condition
on the two-dimensional SSM, we observe in Figure \ref{fig:IC1} that
the full solution indeed stays on the SSM, thereby showing that the
SSM is indeed invariant. Since the computed $\mathbf{W}$ is only
a third-order appoximation of the SSM, an initial condition far enough
from the the fixed point may not stay over this approximate surface
(cf. Figure \ref{fig:IC2}). Indeed, a higher-order approximation
to $\mathbf{W}$ would be required to verify the invariance of $W(\mathcal{E})$
numerically for larger initial conditions. 
\begin{figure}[H]
\subfloat[\label{fig:IC1}Initial condition with polar coordinates $(\rho(0),\theta(0))=(3,0)$]{\centering{}\includegraphics[width=0.49\textwidth]{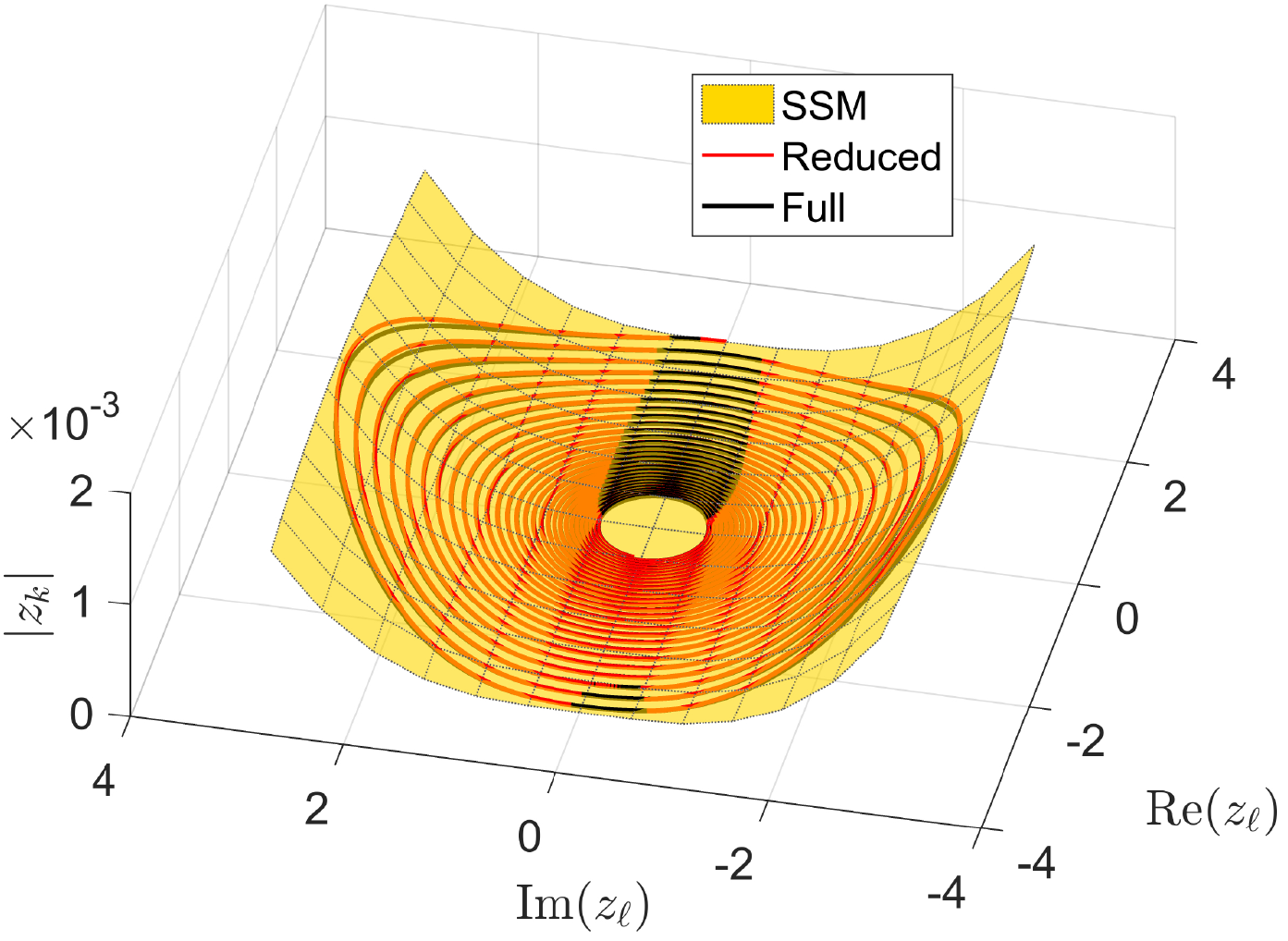}}\hfill{}\subfloat[\label{fig:IC2}Initial condition with polar coordinates $(\rho(0),\theta(0))=(9.0)$]{\centering{}\includegraphics[width=0.49\textwidth]{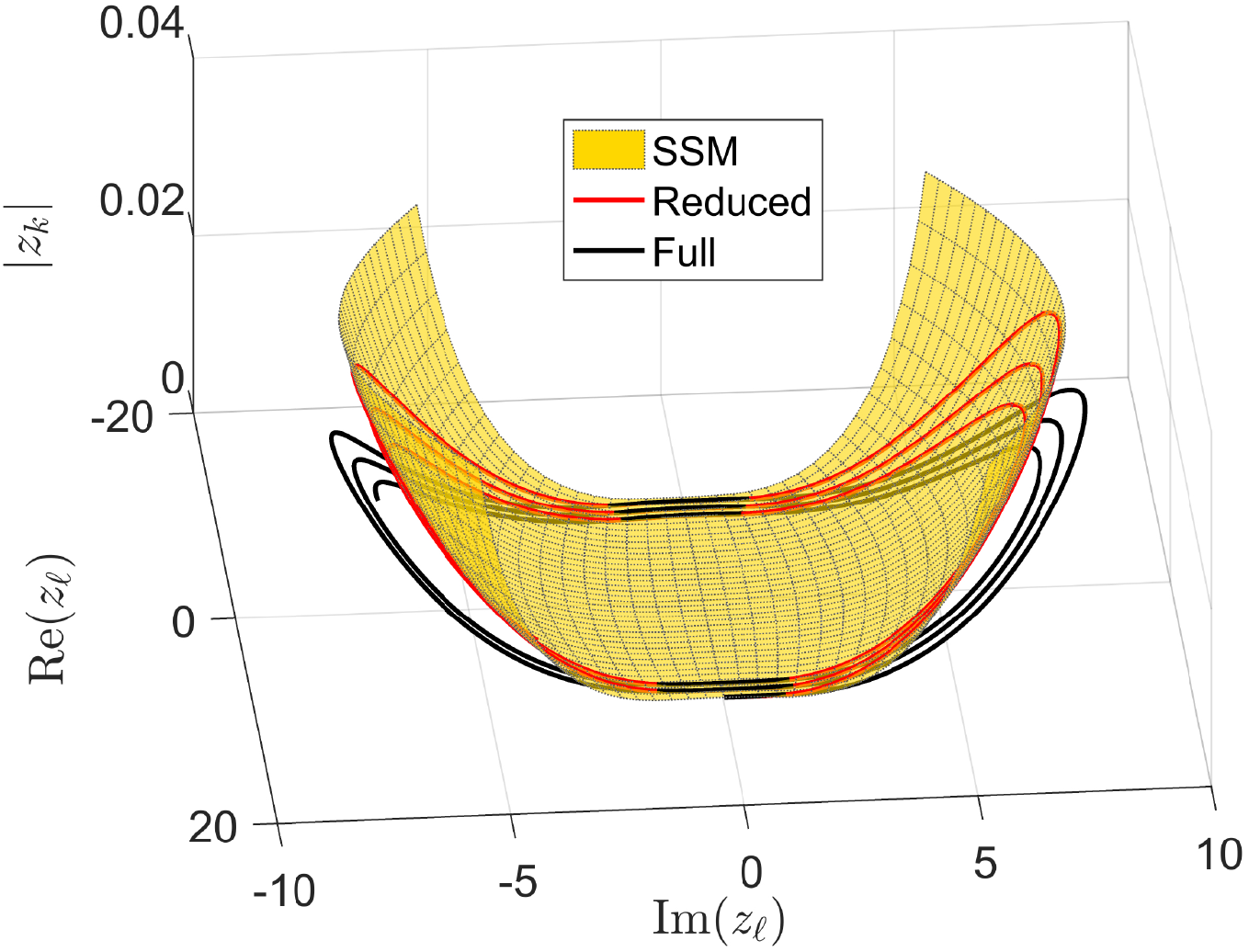}}

\caption{\label{fig:IC-off-SSM-1}The SSM plotted in the modal coordinates
$\mathbf{z}$. Full solution trajectory is shown for an initial condition
taken on the SSM. Modal amplitude $|z_{k}|$ of the $k^{\mathrm{th}}$
mode ($k=3$) is plotted against the master modal variables ($\mathrm{Re}(z_{\ell}),\,\mathrm{Im}(z_{\ell})$)
with $\ell=1$. The full solution trajectory (a) stays approximately
on the computed SSM, (b) does not stay exactly on the computed SSM
for an initialization far enough from the fixed point. }
\end{figure}
\item \textbf{Initial condition off the SSM but still on the slow manifold}:
We illustrate that the SSM is attracting nearby trajectories within
the SFD-based slow manifold by launching a trajectory away from the
SSM but still inside the slow manifold $\mathcal{M}_{\epsilon}(\tau)$.
As shown in Figure \ref{fig:IC-off-SSM}, a trajectory of the full
solution which starts at an arbitrarily chosen point in the phase
space of system (\ref{eq:statespaceROM}), close enough to the computed
SSM, quickly converges towards the SSM and synchronizes with the flow
on the SSM. Figure \ref{fig:exp-decay} shows that this rate of attraction
towards the SSM is indeed faster than the typical decay rate within
the SSM. 
\begin{figure}[H]
\subfloat[]{\centering{}\includegraphics[width=0.49\textwidth]{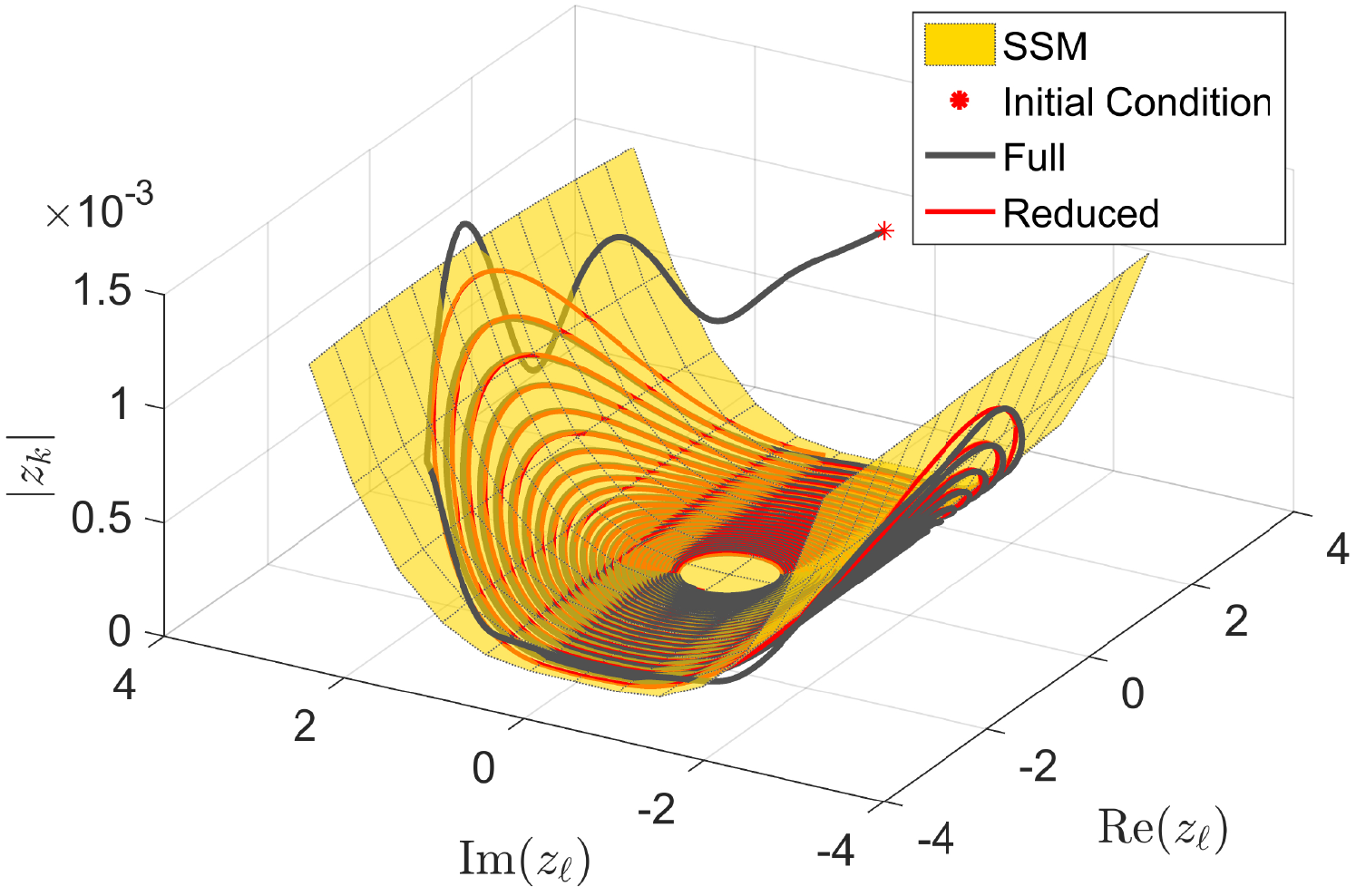}}\hfill{}\subfloat[]{\centering{}\includegraphics[width=0.49\textwidth]{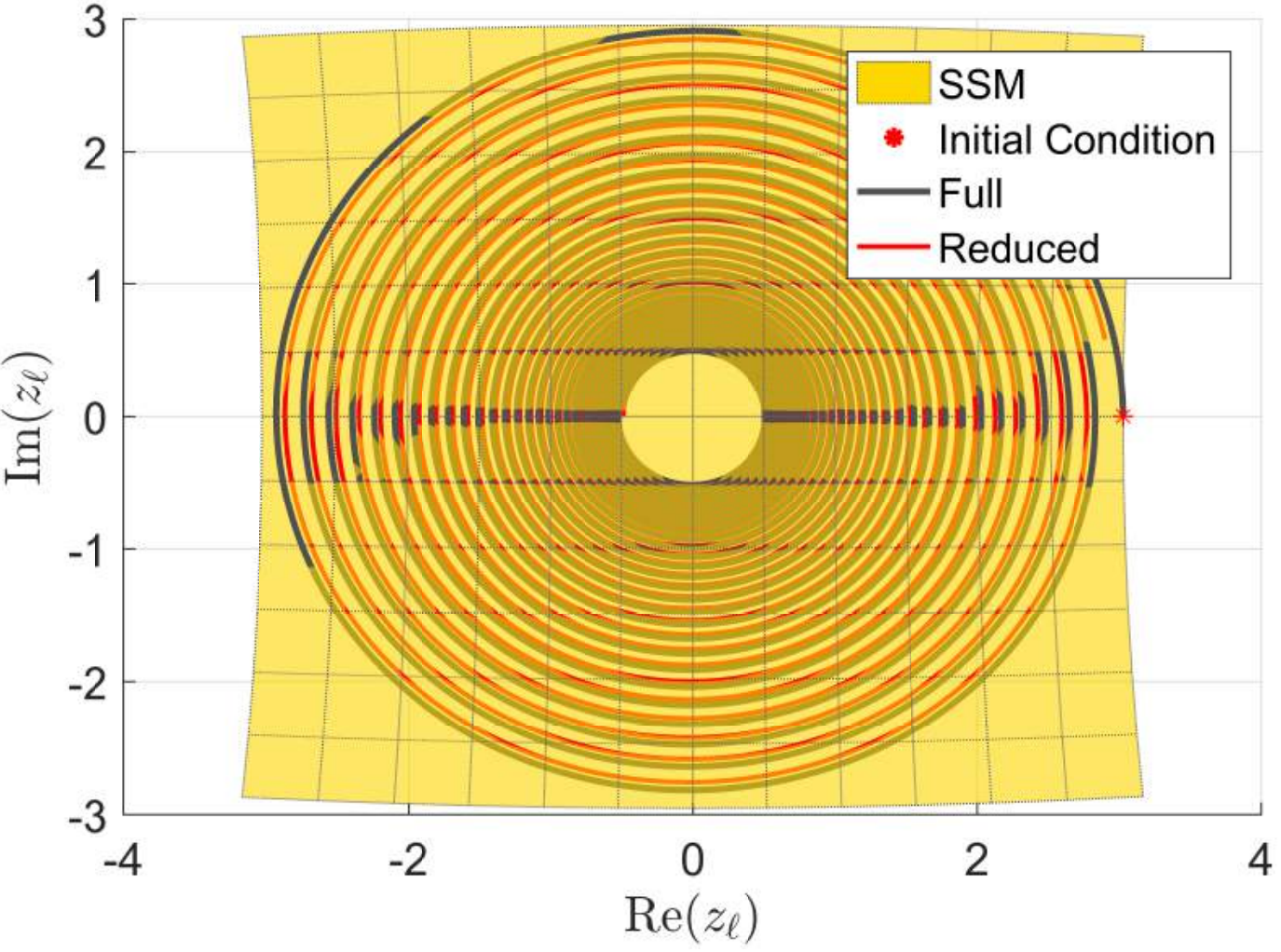}}

\caption{\label{fig:IC-off-SSM} The SSM plotted in the modal coordinates $\mathbf{z}$.
Full solution trajectory is shown for an initial condition taken off
the SSM but still inside the slow manifold. Modal amplitude $|z_{k}|$
of the $k^{\mathrm{th}}$ mode ($k=3$) is plotted against the master
modal unknows ($\mathrm{Re}(z_{\ell}),\,\mathrm{Im}(z_{\ell})$).
The full solution trajectory (a) quickly decays onto the SSM, and
(b) synchronises with the dynamics on the SSM. }
\end{figure}
\begin{figure}[H]
\centering{}\includegraphics[width=0.6\textwidth]{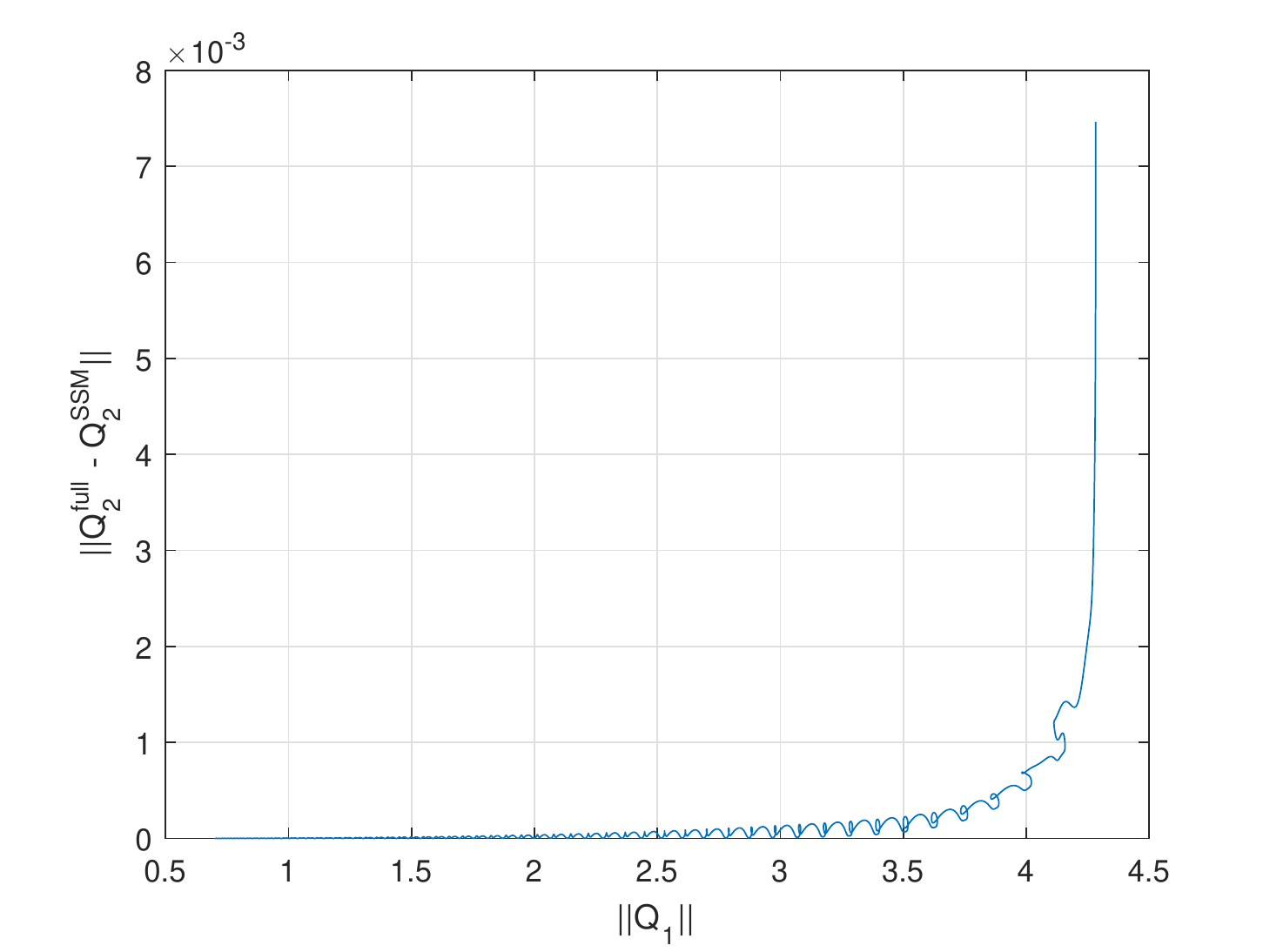}\caption{\label{fig:exp-decay}The exponential decay of trajectories towards
the SSM within the slow manifold in the modal phase space: $|Q_{2}(\mathbf{z})-Q_{2}(\mathbf{W}(\mathbf{s}))|$
vs. $|Q_{1}(\mathbf{z}))|$. The decay rate transverse to the SSM
is exponentially faster than that the typical rates along the SSM.}
\end{figure}
\item \textbf{Initial condition off the slow manifold}: Although the SSM
is computed for the SFD-reduced system (\ref{eq:SFD-ROM-beam}), the
reduced dynamics on the SSM also captures the asymptotic behavior
of the full system by construction. Indeed, as seen in Figure \ref{fig:IC-off-SM},
a trajectory initialized off the slow manifold in the phase space
of system (\ref{eq:FEdisc}) (thus off the SSM as well) is also attracted
towards the SSM, and synchronizes with the corresponding flow on the
SSM. This synchronization can be seen in more detail in Figure \ref{fig:slow-fast-comparison},
where the time histories of slow and fast variables are compared between
the slow and full solution. As expected, the response of the fast
degree of freedom, in particular, shows that the full solution performs
rapid oscilations around, and stabilizes on, the reduced solution
obtained from SFD. Finally, the full solution approaches the SSM-reduced
solution. 
\begin{figure}[h]
\begin{centering}
\includegraphics[width=0.8\textwidth]{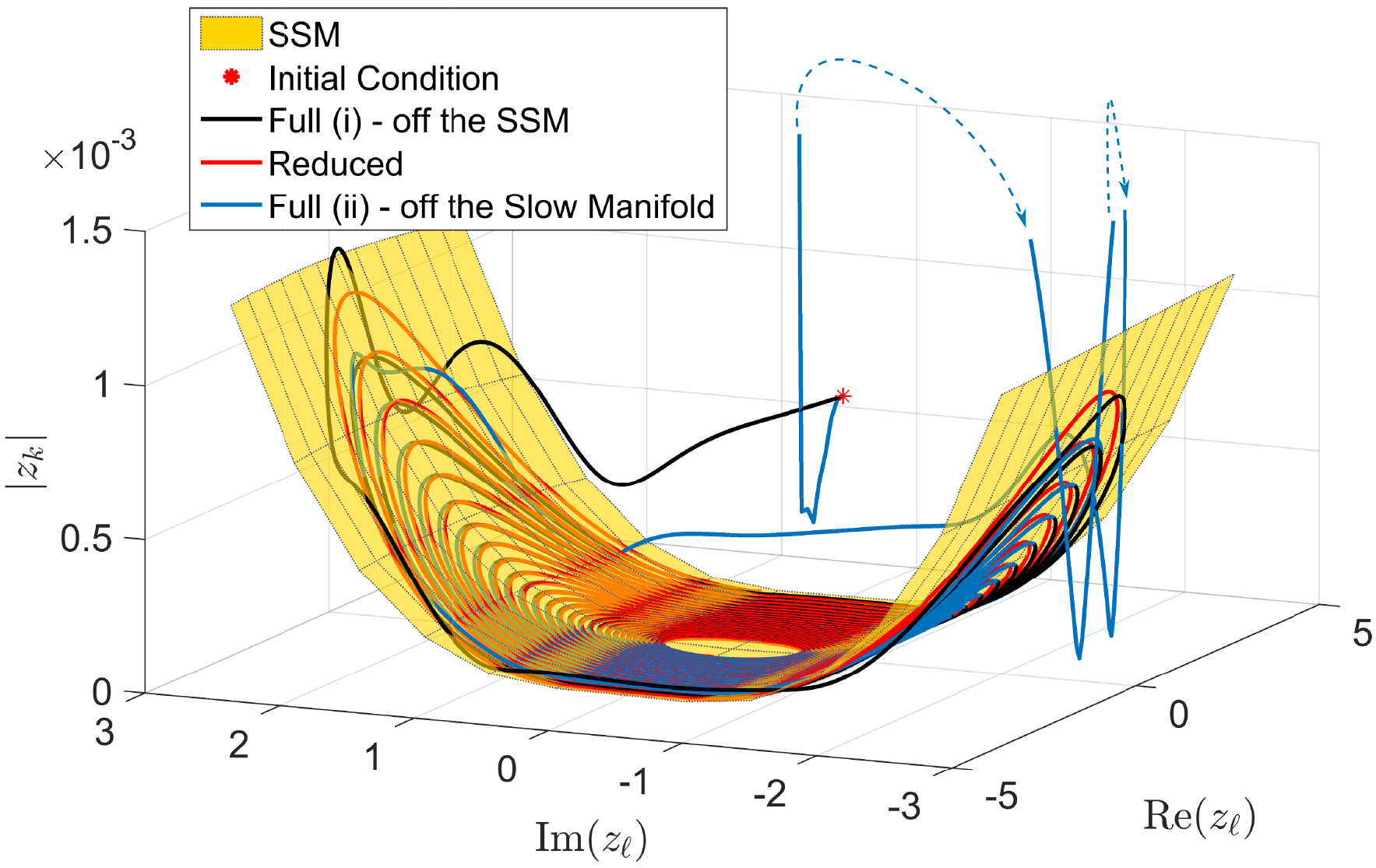}\caption{\label{fig:IC-off-SM}Modal amplitude $|z_{k}|$ of the $k^{\mathrm{th}}$
mode ($k=3$) is plotted against the master modal variables ($\mathrm{Re}(z_{\ell}),\,\mathrm{Im}(z_{\ell})$).
The full (ii) solution trajectory (blue) initialized off the slow
manifold synchronizes with the reduced dynamics on the SSM. The full
(i) soluion trajectory (black) is initialized on the slow manifold
by projecting the original initial confition onto the slow manifold.
Note that in this plot, the two initial conditions appear to coincide,
since the axes feature only the slow (modal) DOFs. The dashed lines
are a schematic representation of the trajectory leaving and returning
to the plot region. }
\par\end{centering}
\end{figure}
\begin{figure}[H]
\subfloat[]{\centering{}\includegraphics[width=0.49\textwidth]{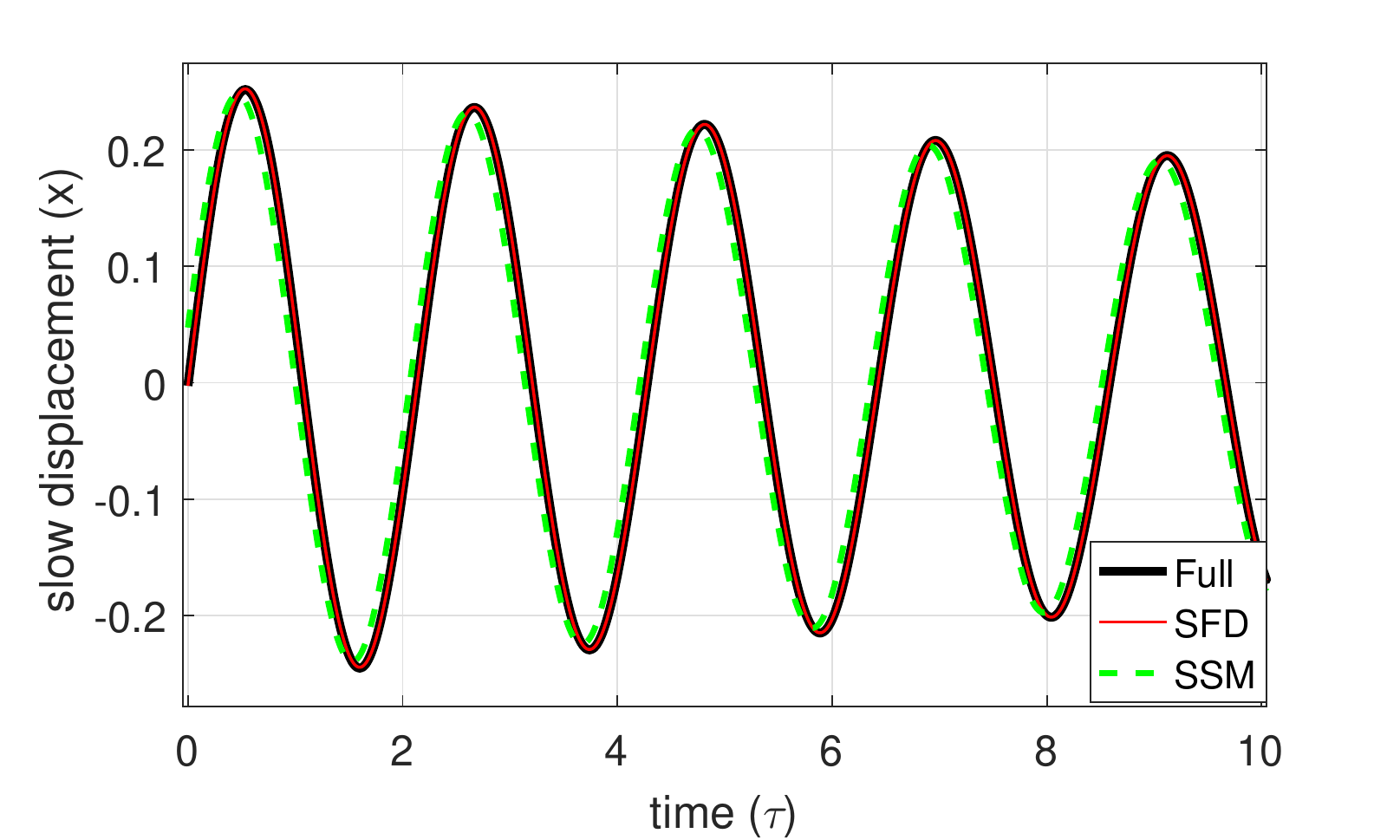}}\hfill{}\subfloat[]{\centering{}\includegraphics[width=0.49\textwidth]{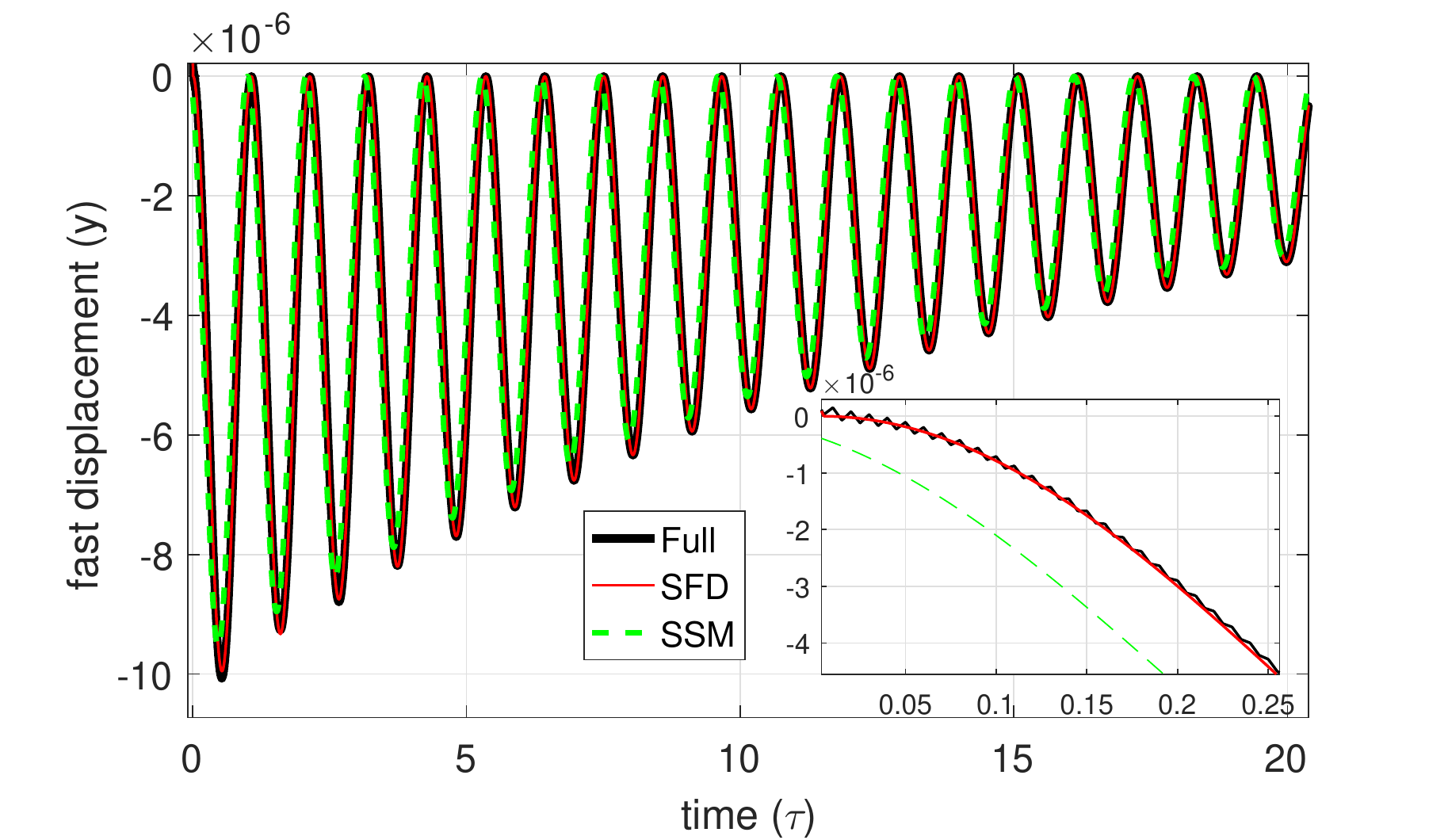}}

\caption{\label{fig:slow-fast-comparison}The comparison of full and different
reduced solutions for (a) slow and (b) fast degrees of freedom. Note
the two timescales in the zoom-in for the dynamics of the fast DOF.
The full solution quickly decays to the slow manifold (with SFD-reduced
solutions) after which it decays to the SSM. }
\end{figure}
\end{enumerate}
Next, we study the rate of decay of the enslaved variables towards
the computed SSM. From Figure \ref{fig:decayRate}, we see that the
initial decay rate of a trajectory towards the SSM is approximately
$e^{-0.4607\tau}$ which is actually dominated by the decay rate suggested
by the second slowest eigenvalue, i.e., $\mathrm{Re}(\lambda_{2})\approx-0.4724$.
The slow dynamics within the SSM is expected to occur at the rate
suggested by the first eigenvalue, i.e., $\mathrm{Re}(\lambda_{1})\approx-0.0295$.
The final decay-rate for the full solution is, however, a decay rate
of $e^{-0.0891\tau}$ which is between the slowest and second slowest
rates, still being an order of magnitude faster than the dynamics
within the SSM, as expected from the underlying theory.

\section{Conclusions}

We have demonstrated two exact model reduction techniques, the slow-fast
decomposition (SFD) and spectral-submanifold (SSM)-based reduction,
on the finite-element model of a von Kármán beam. The SFD enabled
us to express the fast axial variables of the beam as a graph over
the slow transverse variables, yielding a reduced-order model only
in terms of the slow variables. This SFD-reduced model possesses a
gap in its spectral quotient disctribution, indicating an opportunity
for further reduction using the SSM. Subsequently, we carried out
a single-mode SSM-based reduction on the SFD-reduced model. This two-stage,
exact reduction procedure resulted in a drastic reduction of the model
dimension, and ensured that the full system trajectories synchronized
with the reduced model trajectories at rates much faster than typical
decay rates within the manifold.

The application of the above mentioned exact reduction techniques
to our beam model confirms their potential for truly high-dimensional
systems. However, significant work remains to be done to enable application
of these techniques to simulations of realistic structures, especially
on the computational implementation side. In particular, these techniques,
in their present form, require the nonlinear coefficients of the system
to be known apriori. For finite-element-based applications, often
these nonlinear coefficients are embedded in the software at the element
level and their full counterparts are never calculated during a simulation
(for reasons of computational efficiency). This is certainly a computational
challenge for application of SFD/SSM-based reduction techniques to
real-world applications and shall be the focus of our future endeavours.

The current work features a single-mode SSM. This was justified in
our beam example, because decay rates in the first mode were exceptionally
slower than the rest. For more general damped-mechanical systems,
however, we expect higher-dimensional SSMs (i.e., smoothest nonlinear
continuations of high-dimensional linear modal subspaces) to be more
relevant for reduction. Furthermore, the polynomial expansion of the
SSM mapping contains terms only up to the third order in the current
work. As also discussed in Section 4.3, higher-order terms would play
an important role in capturing trajectories initialized farther from
the fixed point. A MATLAB computational toolbox to compute multi-mode
SSMs with terms of up to arbitrary order is currently under development
(cf. Ponsioen and Haller \cite{Ponsioen17}). 

The current work covers the autonomous (non-forced) beam model. Existence
and uniqueness results for SSMs are also available for non-autonomous,
damped mechanical systems as discussed by Haller and Ponsioen \cite{haller16}.
For (quasiperiodically) forced systems, the SSM appears as a quasiperiodically
deforming invariant surface in the phase space, which is again an
ideal tool for model reduction. The computation of such non-autonomous
SSMs will find useful application in non-autonomous mechanical systems,
and is subject to ongoing work (cf. Breunung and Haller \cite{Breunung17}).

\appendix

\section{$\mathcal{O}\left(\epsilon^{2}\right)$terms in the slow manifold}

\label{sec:appA}

As shown by Haller and Ponsioen \cite{haller16sf}, the enslavement
of the fast variables to the slow ones along the slow manifold is
given by the functions

\begin{align*}
\mathbf{y} & =\epsilon\mathbf{G}_{0}(\mathbf{x},\dot{\mathbf{x}},\tau)+\epsilon^{2}\mathbf{G}_{1}(\mathbf{x},\dot{\mathbf{x}},\tau)+\epsilon^{3}\mathbf{G}_{2}(\mathbf{x},\dot{\mathbf{x}},\tau)+\mathcal{O}\left(\epsilon^{4}\right)\,,\\
\dot{\mathbf{y}} & =\epsilon\mathbf{H}_{0}(\mathbf{x},\dot{\mathbf{x}},\tau)+\epsilon^{2}\mathbf{H}_{1}(\mathbf{x},\dot{\mathbf{x}},\tau)+\epsilon^{3}\mathbf{H}_{2}(\mathbf{x},\dot{\mathbf{x}},\tau)+\mathcal{O}\left(\epsilon^{4}\right)\,,
\end{align*}
where $\mathbf{G}_{0},\,\mathbf{G}_{1},\,\mathbf{H}_{0}$ are as shown
in (\ref{eq:OepsTerms}) and $\mathbf{H}_{1}$ is given by 
\begin{equation}
\mathbf{H}_{1}(\mathbf{x},\dot{\mathbf{x}},\tau)=\left[\partial_{\mathbf{x}}\mathbf{G}_{1}\right]\dot{\mathbf{x}}+\left[\partial_{\dot{\mathbf{x}}}\mathbf{G}_{1}\right]\overline{\mathbf{P}_{1}}+\partial_{\tau}\mathbf{G}_{1}+\partial_{\dot{\mathbf{x}}}\mathbf{G}_{0}\left(\overline{\partial_{\boldsymbol{\eta}}\mathbf{P}_{1}}\mathbf{G}_{1}+\overline{\partial_{\dot{\mathbf{y}}}\mathbf{P}_{1}}\mathbf{H}_{0}+\overline{\partial_{\epsilon}\mathbf{P}_{1}}\right).\label{eq:H1}
\end{equation}
In order to obtain the $\mathcal{O}\left(\epsilon^{2}\right)$ terms
in the reduced model, we have calculated the general expressions for
$\mathbf{G}_{2}$ as

\begin{align}
\mathbf{G}_{2}(\mathbf{x},\dot{\mathbf{x}},\tau) & =\left[\overline{\partial_{\boldsymbol{\eta}}\mathbf{P}_{2}}\right]^{-1}\left(\left[\partial_{\mathbf{x}}\mathbf{H}_{0}\right]\dot{\mathbf{x}}+\left[\partial_{\dot{\mathbf{x}}}\mathbf{H}_{0}\right]\overline{\mathbf{P}_{1}}+\partial_{\tau}\mathbf{H}_{0}-\left(\mathbf{J}\mathbf{G}_{1}+\mathbf{L}\mathbf{H}_{0}+\overline{\partial_{\epsilon}^{2}\mathbf{P}_{2}}+\overline{\partial_{\dot{\mathbf{y}}}\mathbf{P}_{2}}\mathbf{H}_{1}\right)\right),\label{eq:G2}
\end{align}
where 
\begin{align*}
\mathbf{J} & =\frac{1}{2}\left(\overline{\partial_{\boldsymbol{\eta}}^{2}\mathbf{P}_{2}}\mathbf{G}_{1}+\overline{\partial_{\dot{\mathbf{y}}}\partial_{\boldsymbol{\eta}}\mathbf{P}_{2}}\mathbf{H}_{0}+\overline{\partial_{\epsilon}\partial_{\boldsymbol{\eta}}\mathbf{P}_{2}}\right)\,,\\
\mathbf{L} & =\frac{1}{2}\left(\overline{\partial_{\dot{\mathbf{y}}}^{2}\mathbf{P}_{2}}\mathbf{H}_{0}+\overline{\partial_{\boldsymbol{\eta}}\partial_{\dot{\mathbf{y}}}\mathbf{P}_{2}}\mathbf{G}_{1}+\overline{\partial_{\epsilon}\partial_{\dot{\mathbf{y}}}\mathbf{P}_{2}}\right)\,.
\end{align*}
The reduced-order model including the $\mathcal{O}\left(\epsilon^{2}\right)$
terms is then given as
\[
\ddot{\mathbf{x}}-\overline{\mathbf{P}_{1}}-\epsilon\left[\overline{\partial_{\boldsymbol{\eta}}\mathbf{P}_{1}}\mathbf{G}_{1}+\overline{\partial_{\dot{\mathbf{y}}}\mathbf{P}_{1}}\mathbf{H_{0}}+\overline{\partial_{\epsilon}\mathbf{P}_{1}}\right]-\epsilon^{2}\left[\mathbf{N}\mathbf{G}_{1}+\overline{\partial_{\boldsymbol{\eta}}\mathbf{P}_{1}}\mathbf{G}_{2}+\mathbf{R}\mathbf{H}_{0}+\overline{\partial_{\dot{\mathbf{y}}}\mathbf{P}_{1}}\mathbf{H}_{1}+\overline{\partial_{\epsilon}^{2}\mathbf{P}_{1}}\right]+\mathcal{O}(\epsilon^{3})=\mathbf{0},
\]
where
\begin{align*}
\mathbf{N} & =\frac{1}{2}\left(\overline{\partial_{\boldsymbol{\eta}}^{2}\mathbf{P}_{1}}\mathbf{G}_{1}+\overline{\partial_{\dot{\mathbf{y}}}\partial_{\boldsymbol{\eta}}\mathbf{P}_{1}}\mathbf{H}_{0}+\overline{\partial_{\epsilon}\partial_{\boldsymbol{\eta}}\mathbf{P}_{1}}\right)\,,\\
\mathbf{R} & =\frac{1}{2}\left(\overline{\partial_{\dot{\mathbf{y}}}^{2}\mathbf{P}_{1}}\mathbf{H}_{0}+\overline{\partial_{\boldsymbol{\eta}}\partial_{\dot{\mathbf{y}}}\mathbf{P}_{1}}\mathbf{G}_{1}+\overline{\partial_{\epsilon}\partial_{\dot{\mathbf{y}}}\mathbf{P}_{1}}\right)\,.
\end{align*}

In our beam framework, the general expressions in (\ref{eq:H1},\ref{eq:G2})
become
\begin{align}
\mathbf{H}_{1}(\mathbf{x},\dot{\mathbf{x}},\tau) & =\mathbf{K}_{2}^{-1}\left(\zeta\left[\partial_{\mathbf{x}}\boldsymbol{\mathcal{H}}(\mathbf{x})-\boldsymbol{\mathcal{E}}(\mathbf{x})\right]\overline{\mathbf{P}_{1}}+\zeta\left[\partial_{\mathbf{x}}^{2}\boldsymbol{\mathcal{H}}(\mathbf{x})-\partial_{\mathbf{x}}\boldsymbol{\mathcal{E}}(\mathbf{x})\right]:(\dot{\mathbf{x}}\otimes\dot{\mathbf{x}})+\beta\dot{\mathbf{p}}\right)\,,\nonumber \\
\mathbf{G}_{2}(\mathbf{x},\dot{\mathbf{x}},\tau) & =\left[\mathbf{K}_{2}\mathbf{M}_{2}^{-1}\mathbf{K}_{2}\right]^{-1}\left(\left[\partial_{\mathbf{x}}^{2}\boldsymbol{\mathcal{H}}(\mathbf{x})\right]:(\dot{\mathbf{x}}\otimes\dot{\mathbf{x}})+\left[\partial_{\mathbf{x}}\boldsymbol{\mathcal{H}}(\mathbf{x})\right]\overline{\mathbf{P}_{1}}\right)-\zeta\mathbf{H}_{1}.\label{eq:H1G2-beam}
\end{align}
For the assumed visco-elastic material damping, the expressions for
$\mathbf{G}_{1}(\mathbf{x},\dot{\mathbf{x}},\tau)$ in (\ref{eq:OepsTerms})
and $\mathbf{H}_{1}(\mathbf{x},\dot{\mathbf{x}},\tau)$ in (\ref{eq:H1G2-beam})
can be further simplified to
\begin{align*}
\mathbf{G}_{1}(\mathbf{x},\dot{\mathbf{x}},\tau) & =\beta\mathbf{K}_{2}^{-1}\mathbf{p}(\tau)\,,\\
\mathbf{H}_{1}(\mathbf{x},\dot{\mathbf{x}},\tau) & =\beta\mathbf{K}_{2}^{-1}\dot{\mathbf{p}}(\tau)\,,
\end{align*}
and the reduced-order model can be simplified as 
\begin{gather}
\mathbf{M}_{1}\ddot{\mathbf{x}}+\mathbf{K}_{1}\mathbf{x}+\boldsymbol{\mathcal{F}}\left(\mathbf{x},\mathbf{G}_{0}(\mathbf{x},\dot{\mathbf{x}},\tau)\right)+\boldsymbol{\mathcal{G}}(\mathbf{x})+\nonumber \\
\epsilon\left[\overline{\partial_{\boldsymbol{\eta}}\boldsymbol{\mathcal{F}}\left(\mathbf{x},\boldsymbol{\eta}\right)}\mathbf{G}_{1}(\mathbf{x},\dot{\mathbf{x}},\tau)+\zeta\left(\boldsymbol{\mathcal{D}}(\mathbf{x})\mathbf{H}_{0}(\mathbf{x},\dot{\mathbf{x}},\tau)+\left(\mathbf{K}_{1}+\boldsymbol{\mathcal{C}}(\mathbf{x})\right)\dot{\mathbf{x}}\right)\right]+\label{eq:SFD-ROM-beam-2}\\
\epsilon^{2}\left[\overline{\partial_{\boldsymbol{\eta}}\boldsymbol{\mathcal{F}}\left(\mathbf{x},\boldsymbol{\eta}\right)}\mathbf{G}_{2}(\mathbf{x},\dot{\mathbf{x}},\tau)+\zeta\boldsymbol{\mathcal{D}}(\mathbf{x})\mathbf{H}_{1}(\mathbf{x},\dot{\mathbf{x}},\tau)\right]+\mathcal{O}\left(\epsilon^{3}\right)=\alpha\mathbf{q}(\tau).\nonumber 
\end{gather}

\section{Forumulae for single-mode SSMs}

\label{sec:appB}

We express the general formulae of Szalai et al. \cite{Szalai16}
for Taylor coefficients of the SSM mapping $\mathbf{W}(\mathbf{s})$
in our present notation. Let $\boldsymbol{\mathcal{T}}(\mathbf{z})$
from (\ref{eq:diagROM}) be a general polynomial of the form 
\begin{equation}
\boldsymbol{\mathcal{T}}(\mathbf{z})=\mathbf{\boldsymbol{\mathcal{T}}}^{(2)}(\mathbf{z})+\mathbf{\boldsymbol{\mathcal{T}}}^{(3)}(\mathbf{z})+\dots,\label{eq:Texpansion2}
\end{equation}
where $\mathbf{\boldsymbol{\mathcal{T}}}^{(n)}(\mathbf{z})$ denotes
the $n^{\mathrm{th}}$order terms of $\boldsymbol{\mathcal{T}}$.
The first order terms in the expansion (\ref{eq:Wexpansion}) are
given by
\[
\mathbf{\boldsymbol{\mathcal{W}}}^{(1)}(\mathbf{s})=\mathbf{W}^{(1)}\mathbf{s}\,,
\]
where $\mathbf{W}^{(1)}\in\mathbb{C}^{2n_{s}\times2}$ is an all-zero
matrix except for two non-zero entries given by $\left(\mathbf{W}^{(1)}\right)_{1,\ell}=\lambda_{\ell},\,\left(\mathbf{W}^{(1)}\right)_{2,\ell+1}=\lambda_{\ell+1}$.
The $i^{\mathrm{th}}$ component of the quadratic terms can be written
as
\[
\left(\mathbf{\boldsymbol{\mathcal{W}}}^{(2)}(\mathbf{s})\right)_{i}=\sum_{j=1}^{2}\sum_{k=1}^{2}W_{ijk}^{(2)}s_{j}s_{k}\,,\quad i\in\{1,\dots,2n_{s}\},\quad j,k\in\{1,2\}\,,
\]
where $\mathbf{W}^{(2)}\in\mathbb{C}^{2n_{s}\times2\times2}$ is a
sparse 3-tensor with nonzero entries given as 
\begin{align*}
W_{i11}^{(2)} & =\frac{T_{i\ell\ell}^{(2)}}{2\lambda_{\ell}-\lambda_{i}}\,,\quad i\in\{1,\dots,2n_{s}\}\,,\\
W_{i22}^{(2)} & =\frac{T_{i(\ell+1)(\ell+1)}^{(2)}}{2\lambda_{\ell+1}-\lambda_{i}}\,,\quad i\in\{1,\dots,2n_{s}\}\,,\\
W_{ijk}^{(2)} & =\frac{T_{ijk}^{(2)}}{\lambda_{\ell}+\bar{\lambda}_{\ell}-\lambda_{i}}\,,\quad i\in\{1,\dots,2n_{s}\},\quad(j,k)\in\{(\ell,\ell+1),(\ell+1,\ell)\}\,.
\end{align*}
Finally, the cubic terms can be written as 
\[
\left(\mathbf{\boldsymbol{\mathcal{W}}}^{(3)}(\mathbf{s})\right)_{i}=\sum_{j=1}^{2}\sum_{k=1}^{2}\sum_{l=1}^{2}W_{ijkl}^{(3)}s_{j}s_{k}s_{l}\,,\quad i\in\{1,\dots,2n_{s}\},\quad j,k,l\in\{1,2\}\,,
\]
where $\mathbf{W}^{(3)}\in\mathbb{C}^{2n_{s}\times2\times2\times2}$
is a sparse 4-tensor with nonzero entries given as
\begin{align*}
W_{i111}^{(3)} & =\frac{\sum_{j=1}^{2n_{s}}\left[\left(1+\delta_{\ell j}\right)\left(T_{ij\ell}^{(2)}+T_{i\ell j}^{(2)}\right)W_{j11}^{(2)}\right]+T_{i\ell\ell\ell}^{(3)}}{3\lambda_{\ell}-\lambda_{i}}\,,\quad i\in\{1,\dots,2n_{s}\}\,,\\
W_{i222}^{(3)} & =\frac{\sum_{j=1}^{2n_{s}}\left[\left(1+\delta_{(\ell+1)j}\right)\left(T_{ij(\ell+1)}^{(2)}+T_{i(\ell+1)j}^{(2)}\right)W_{j22}^{(2)}\right]+T_{i(\ell+1)(\ell+1)(\ell+1)}^{(3)}}{3\lambda_{\ell+1}-\lambda_{i}}\,,\quad i\in\{1,\dots,2n_{s}\}\,,\\
W_{ijkl}^{(3)} & =(1-\delta_{i\ell})\frac{V_{i}+T_{ijkl}^{(3)}}{2\lambda_{\ell}+\bar{\lambda}_{\ell}-\lambda_{i}}\,,\quad i\in\{1,\dots,2n_{s}\},\quad(j,k,l)\in\{(\ell,\ell+1,\ell),(\ell,\ell,\ell+1),(\ell+1,\ell,\ell)\}\,,\\
W_{ijkl}^{(3)} & =(1-\delta_{i(\ell+1)})\frac{U_{i}+T_{ijkl}^{(3)}}{2\lambda_{\ell+1}+\bar{\lambda}_{\ell+1}-\lambda_{i}}\,\\
 & \quad i\in\{1,\dots,2n_{s}\},\quad(j,k,l)\in\{(\ell+1,\ell,\ell+1),(\ell+1,\ell+1,\ell),(\ell,\ell+1,\ell+1)\}\,,
\end{align*}
with 
\begin{align*}
V_{i} & =\sum_{j=1}^{2n_{s}}\left[\left(1+\delta_{\ell j}\right)\left(T_{ij\ell}^{(2)}+T_{i\ell j}^{(2)}\right)\left(W_{j21}^{(2)}+W_{j12}^{(2)}\right)+\left(1+\delta_{(\ell+1)j}\right)\left(T_{ij(\ell+1)}^{(2)}+T_{i(\ell+1)j}^{(2)}\right)W_{j11}^{(2)}\right]\,,\quad i\in\{1,\dots,2n_{s}\},\\
U_{i} & =\sum_{j=1}^{2n_{s}}\left[\left(1+\delta_{(\ell+1)j}\right)\left(T_{ij(\ell+1)}^{(2)}+T_{i(\ell+1)j}^{(2)}\right)\left(W_{j21}^{(2)}+W_{j12}^{(2)}\right)+\left(1+\delta_{\ell j}\right)\left(T_{ij\ell}^{(2)}+T_{i\ell j}^{(2)}\right)W_{j22}^{(2)}\right]\,,\quad i\in\{1,\dots,2n_{s}\}.
\end{align*}
Here $T_{ijkl}^{(3)},T_{ijk}^{(2)}\in\mathbb{C}$ denote the components
of the 4-tensor and 3-tensor, respectively, in the expansion of $\boldsymbol{\mathcal{T}}$
given in (\ref{eq:Texpansion2}), and $\delta_{ij}$ respresents the
Kronecker-delta. Furthermore, the general expression for $\beta_{\ell}$
used in the expansion for $\mathbf{R}$ as in (\ref{eq:SSMROM}),
as derived by Szalai et al. \cite{Szalai16}, can be written in our
notation as
\[
\beta_{\ell}=V_{\ell}+T_{\ell\ell\ell(\ell+1)}+T_{\ell\ell(\ell+1)\ell}+T_{\ell(\ell+1)\ell\ell}.
\]
Note that the Einstein summation convention has \textit{not} been
followed in any of the above expressions.

\end{document}